\documentclass[sort&compress]{elsarticle}
 \pdfoutput=1 
 \journal{Journal of Computational Physics}
\usepackage{amsmath,amssymb,graphicx,mathrsfs,subfigure}
\usepackage{amsmath,mathtools}
\usepackage[normalem]{ulem}
\usepackage[usenames,dvipsnames]{xcolor}
 \definecolor{burgundy}{rgb}{0.59, 0.0, 0.09}

\usepackage{hyperref}
\hypersetup{hidelinks}

\newcommand{\bm}[1]{\mbox{\boldmath $ #1 $}}    
\renewcommand{\vec}[1]{\bm{#1}}                 
\newcommand{\mat}[1]{\mathbf{#1}}		
\newcommand{\ten}[1]{\mathcal{#1}}		
\DeclareMathOperator{\Tr}{\textup{tr}}		
\DeclareMathOperator{\Det}{\textup{det}}		
\newcommand{\norm}[1]{\left\lVert #1 \right\rVert} 
\newcommand{\fraco}[2]{\displaystyle{\frac{#1}{#2}}}

\newcommand{\C}{\ten{C}}
\newcommand{\G}{\ten{G}}
\newcommand{\Gz}{\ten{G}_0}

\newcommand{\lf}{\left}
\newcommand{\rt}{\right}

\def\onedot{$\mathsurround0pt\ldotp$}
\def\cddot{
  \mathbin{\vcenter{\baselineskip.67ex
    \hbox{\onedot}\hbox{\onedot}}%
  }}%
\begin{document}
\begin{frontmatter}
\title{A continuous energy-based immersed boundary method for elastic shells}
\author[om,ws]{Ondrej Maxian}
\ead{oxm60@case.edu}
\author[ak]{Andrew T. Kassen}
\ead{kassen@math.utah.edu}
\author[ws]{Wanda Strychalski\corref{cor}}
\ead{wis6@case.edu}
\cortext[cor]{Corresponding Author}

\address[om]{Department of Chemical and Biomolecular Engineering,
Case Western Reserve University, Cleveland, OH  44106}
\address[ws]{Department of Mathematics, Applied Mathematics, and Statistics, 
Case Western Reserve University, Cleveland, OH  44106}
\address[ak]{Department of Mathematics, University of Utah, UT  84112}

%
%
\begin{abstract}
%
%
The immersed boundary method is a mathematical formulation and numerical method for solving fluid-structure interaction problems. For many biological problems, such as models that include the cell membrane, the immersed structure is a two-dimensional infinitely thin elastic shell immersed in an incompressible viscous fluid. When the shell is modeled as a hyperelastic material, forces can be computed by taking the variational derivative of an energy density functional. A new method for computing a continuous force function on the entire surface of the shell is presented here. The new method is compared to a previous formulation where the surface and energy functional are discretized before forces are computed. 
For the case of Stokes flow, a method for computing quadrature weights is provided to ensure the integral of the elastic spread force density remains zero throughout a dynamic simulation. 
Tests on the method are conducted and show that it yields more accurate force computations than previous formulations as well as more accurate geometric information such as mean curvature. The method is then applied to a model of a red blood cell in capillary flow and a 3D model of cellular blebbing. 
\end{abstract}


\end{frontmatter}
%
%
%
\section{Introduction}
%
%
The immersed boundary (IB) method was first introduced by Charles Peskin in 1972 to model blood flow around heart valves \cite{peskin1972flow}. Since then, it has been applied to a wide variety of biological models, including insect flight \cite{miller2005computational,miller2004vortices}, animal swimming \cite{fauci1988computational}, and cellular mechanics \cite{strychalski2015poroelastic}. 
In  three-dimensional models, computing the elasticity of an immersed 2D elastic membrane or 3D elastic solid as it deforms is more challenging  than in 2D because it involves stress tensors. The first studies in this area used lattice-spring models of elasticity (e.g. \cite{beyer1992computational}), but these methods are disconnected from a constitutive law. The finite element method has also been used to compute the elasticity of immersed structures 
 \cite{griffith2016hybrid, boffi2003finite, boffi2008hyper,boffi2007cfl}, but has not been extended to infinitely thin shells. The IB method has been applied to hyperelastic solids or shells whose energies are governed by a strain energy density functional \cite{deven,fai}.  
However, the surface or solid representation has been assumed to be piecewise linear, and the accuracy of force computations with such a method has not been rigorously tested.

The authors showed in \cite{deven} that for a hyperelastic solid, forces can be computed without using stress tensors. In this work, the derivation was limited to 
Cartesian coordinates and applied to solids and thick elastic shells. The method was later applied to infinitely thin elastic shells represented in two curvilinear coordinates and was 
subsequently used to model red blood cells \cite{fai}, osmotic swelling \cite{wu2015simulation}, and two-phase gels \cite{camacho2016eulerian}. However, these models have used a piecewise linear boundary representation (surface triangulation) \cite{fai}  or marker and cell discretization of the surface \cite{camacho2016eulerian}. Both approaches are based on the assumption that curved surfaces are locally planar, which introduces surface discretization error and fails to capture the inherent curvature of the surface.

Studies on surface representations have generally been conducted separately from those on force computations. In \cite{shankar,shankar2015augmenting}, the authors focus on representing thin surfaces continuously with spherical harmonic or radial basis function (RBF) interpolants. In both cases,  force computations are treated through simple explicit expressions for surface tension \cite{shankar} or fiber elasticity \cite{shankar2015augmenting}. The authors do not consider finding forces through variational derivatives.

Boundary integral methods (BIMs) have been progressing concurrently with the IB method and have been used by others to model fluid-structure interactions in zero Reynolds number flow \cite{kumar2015cell, french,pozrikidis1992boundary, pozrikidis2001interfacial}. In BI methods, the ``hydrodynamic traction jump" across the membrane \cite{kumar2015cell}, which is equivalent to the Lagrangian elastic force density per unit current area \cite{pozrikidis2001interfacial}, must be computed independently prior to integration.
\cite{zhao2010spectral} and \cite{veerapaneni2011fast} both demonstrate the use of a spectrally convergent spherical harmonic representation for surfaces that allows for the calculation of force in BIMs. However, the computations in \cite{zhao2010spectral} and \cite{veerapaneni2011fast} both follow the traditional formulation of BI methods in computing the traction jump from a system of equations involving the Cauchy stress tensor, which itself comes from a series of tensor expansions. Additionally, BI methods require numerical schemes to resolve integrals that are singular at the immersed surfaces. 

The aim of this work is to provide a ``bridge" between the fields of surface representation and energy-based force functions within the IB method. Our goal is to use the continuous surface representations from \cite{shankar} to formulate a continuous function that represents the force density on the membrane as a function of its curvilinear coordinates. In this manner, we avoid the discretization error from surface representation. The advantages of the approach presented here are that the force is computed directly in a single computation, bypassing the need to compute stress tensors at each time step, the interpolation of structure velocities in the IB method occurs over a coarser mesh, and the use of the IB method avoids the need for numerical methods to compute singular integrals as in BI methods. 

Since many biological phenomena  occur on the cellular level, the small length scales (microns) lead to very small Reynolds numbers (Re $\sim 10^{-5}-10^{-2}$). In order to simulate Stokes equations with periodic boundary conditions, the integral of the force density over the fluid domain must be zero \cite{batchelor1967introduction}. In \cite{teran2009tether}, the authors proposed a method for ensuring this condition is met when tether forces are used in the IB method. For translation invariant hyperelastic materials represented by continuous functions, the integral condition is satisfied in the continuous formulation \cite{peskinActa}. However, extra care must be taken to ensure that the discrete integral of the forces in the IB method is zero. Another feature of this work is that we provide a method for computing the correct quadrature weights for closed surfaces so that the IB method can be used with 3D Stokes flow.

The rest of the paper is organized as follows. In Section \ref{sec:mathform}, we present the mathematical formulation of the immersed membrane problem, including our method of computing force densities and correct area quadrature weights to ensure that the forces on the fluid integrate to zero in Stokes' flow. We continue by describing an alternative method based on piecewise linear surface discretization, which we use for comparison. In Section \ref{sec:numform}, we discuss the numerical formulation and discretization of the problem and provide implementation details. We conclude by presenting our results in Section \ref{sec:results}, where we perform force computations and IB simulations on two test objects. In Section \ref{sec:apps}, we apply our model to two different biological processes: a red blood cell flowing through a capillary and a 3D model of cellular blebbing.

%
%
%
\section{Mathematical Formulation}
\label{sec:mathform}
In this section, we begin by presenting an overview of the model equations and mathematical formulation. We then  describe how a continuous elastic force density function can be computed on shells using the energy-based formulation. We provide a method for computing quadrature weights so that the forces applied to the fluid integrate to zero. We compare this to the discretized surface formulation from \cite{kraus1996fluid, fai} and conclude by presenting the energy functionals used in our numerical tests. 
\subsection{Immersed boundary method formulation}
%
%
Our model system is an infinitely thin elastic shell immersed in a viscous fluid. For our applications, the shell represents the cell membrane. Due to the small length scales on the cellular level, inertial forces can be neglected and the fluid obeys Stokes' equations, 
\begin{gather}
\label{eq:Stokeseqnsone}
\mu \Delta \vec{u} - \nabla p + \vec{f} = \vec{0}, \\
\label{eq:Stokeseqnstwo}
\nabla \cdot \vec{u} = 0,
\end{gather}
where $\vec{u}$ represents the fluid velocity, $p$ is the pressure, $\vec{f}$ is an external force, and $\mu$ is the fluid viscosity. The model is formulated using the IB method so that the thin elastic shell is represented on a moving Lagrangian coordinate system while the velocity and pressure are represented on a fixed Eulerian grid. Elastic forces are computed on the Lagrangian grid and spread onto the Eulerian grid to construct the external force density. The spreading operator S is defined as
\begin{equation}
\label{eq:spread}
\vec{f} = \textrm{S} \vec{F}  = \int_{\Gamma}  \vec{F}(\vec{q},t) \delta (\vec{x} - \vec{X}(\vec{q},t))\, d \vec{q},
\end{equation}
where $\Gamma$ represents the Lagrangian structure (shell), $\vec{q} = (q_1,q_2)$ represents material surface coordinates, $\vec{X}(\vec{q},t)$ represents the position of the  structure at time $t$, and $\delta(\vec{x})$ represents the Dirac delta function. We use the convention that lower case letters indicate quantities on the Eulerian grid and capitalized letters denote quantities on the Lagrangian structure.
The structure is updated with the local fluid velocity, satisfying the no-slip boundary condition,
\begin{equation}
\frac{d \vec{X}}{dt} = \vec{U},
\end{equation}
where $\vec{U}$ is the interpolated fluid velocity. The interpolation operator is the adjoint of the spreading operator in 
Eq.\ \eqref{eq:spread} and is given by
\begin{equation}
\vec{U} =\textrm{S}^* \vec{u} = \int_{\Omega} \vec{u} (\vec{x},t) \delta(\vec{x}-\vec{X}(\vec{q},t) )\, d\vec{x},
\end{equation}
where $\Omega$ denotes the fluid domain.
\subsection{Computing the elastic forces}
%
%
We outline two different methods to compute Lagrangian forces analytically based on the assumption that the structure is a hyperelastic material characterized by an energy density functional 
\begin{equation}
\label{eq:energyfunc}
\mathcal{E} = \int_{\Gamma_0} W \, d \vec{q},
\end{equation}
where $W$ is a strain energy density specified by a constitutive law, such as the neo-Hookean model. Exact forms of $W$ are provided in Section \ref{sec:energyfunctionals}. The energy-based model for computing Lagrangian forces on an immersed surface was derived in \cite{deven}. Let the surface $\Gamma$ with undeformed configuration $\Gamma_0$ be defined by Lagrangian coordinates $q_1$ and $q_2$ ($q_i$'s are the parameters). Then the force density per unit reference configuration is related to the variational derivative of energy $\fraco{\delta{\mathcal{E}}}{\delta \vec{X}}$ by
\begin{equation}
\label{eq:forcefromvar}
\vec{F} = -\frac{\delta{\mathcal{E}}}{\delta \vec{X}}(q_1,q_2,t),
\end{equation}
where $\mathcal{E}$ is the surface energy and $t$ is time.

The variational derivative is defined as follows. Let $\phi_{\delta \vec{X}}(\alpha) = \mathcal{E}(\vec{X} + \alpha \delta \vec{X})$. Then the variational derivative is the function $\fraco{\delta{\mathcal{E}}}{\delta \vec{X}}$ that obeys
\begin{equation}
\label{eq:vardef}
\displaystyle{\frac{d \phi_{\delta \bm{X}}(\alpha)}{d\alpha}}\Bigg|_{\alpha=0} = \int_{\Gamma_0} \frac{\delta{\mathcal{E}}}{\delta \vec{X}}(q_1,q_2,t) \cdot \delta \vec{X}(q_1,q_2) \, dA_0.
\end{equation}
In \cite{deven}, this theory was applied to discrete solids in Cartesian coordinates. We extend the theory further to curvilinear coordinates on \textit{continuous closed} surfaces.

In the subsequent sections, we describe two different approaches for computing the force density in Eq.\ \eqref{eq:forcefromvar}. We first outline a procedure to find the variational derivative of the continuous energy density functional and follow with a discussion on computing the proper area weights to obtain force from force density. We next discuss the second method, the linear discrete surface method (LDSM) from \cite{fai}. The difference between the two methods lies primarily in the fact that our formulation relies on a \textit{continuous} force function whose integral is discretized after forces are computed, whereas in LDSM, the surface and energy functional are first discretized. Forces are then computed after the discretized surface energy is exactly integrated.
\subsection{Computing forces from the variational derivative of the energy functional}
\label{sec:vdderivation}
\label{sec:forcevariations}
%
%
Since the surfaces in our applications are topological spheres, we parameterize our surface in terms of spherical coordinates. However, our derivations can be generalized to any coordinate system. 
We begin by defining a two dimensional surface in terms of spatial coordinates $q_1=\lambda$ and $q_2=\theta$. Let
  \begin{align}
  \label{eq:surfacereps}
     \vec{Z}(\lambda,\theta) =\begin{pmatrix}
         \hat{X}(\lambda,\theta) \\
         \hat{Y}(\lambda,\theta) \\
         \hat{Z}(\lambda,\theta)
        \end{pmatrix} 
        \textrm{  and  }
     \vec{X}(\lambda,\theta) =\begin{pmatrix}
         X(\lambda,\theta) \\
         Y(\lambda,\theta) \\
         Z(\lambda,\theta)
        \end{pmatrix}
  \end{align}
be the reference and current configurations of the surface, respectively, where $-\pi < \lambda \leq \pi$ and  $-\pi/2 < \theta \leq \pi/2$. Here
$\vec{X}, \vec{Z}:( \lambda, \theta) \mapsto \mathbb{R}^3$, and we have intentionally expressed Eq.\ \eqref{eq:surfacereps} in a general form to emphasize that our derivation can be applied to any 
particular choice of $\vec{X}$ that describes a closed, smooth surface. 
%
Define the matrices
  \begin{align}
     \mat{\nabla_{\theta}X}& =\begin{pmatrix}
        \fraco{\partial\bm{X}}{\partial \lambda} & \fraco{\partial\bm{X}}{\partial \theta}
        \end{pmatrix}
        \end{align}
        and
        \begin{align}
        \mat{\nabla_{\theta}Z}& =\begin{pmatrix}
        \fraco{\partial\bm{Z}}{\partial \lambda} & \fraco{\partial\bm{Z}}{\partial \theta}
        \end{pmatrix},
  \end{align}
 where $\fraco{\partial\vec{X}}{\partial \lambda}$ is a vector in $\mathbb{R}^3$ that represents the component-wise partial $\lambda$ derivative of $\vec{X}$ (the other components of $ \mat{\nabla_{\theta}X}$ and $ \mat{\nabla_{\theta}Z}$ are defined similarly). The matrix $ \mat{\nabla_{\theta}X}$ is the deformation gradient and the rectangular analogue of the square matrix $\mathbb{F}$ in \cite{deven}. We define the metric tensors
   \begin{align}
   \label{eq:Gdef}
     \G = (\mathbf{\nabla_{\theta}X})^\textrm{T} (\mathbf{\nabla_{\theta}X}) & =      \begin{pmatrix}
\fraco{\partial\bm{X}}{\partial \lambda} \cdot \fraco{\partial\bm{X}}{\partial \lambda} & \fraco{\partial\bm{X}}{\partial \lambda} \cdot \fraco{\partial\bm{X}}{\partial \theta}\\[10 pt]
\fraco{\partial\bm{X}}{\partial \theta} \cdot \fraco{\partial\bm{X}}{\partial \lambda} & \fraco{\partial\bm{X}}{\partial \theta} \cdot \fraco{\partial\bm{X}}{\partial \theta}\\
        \end{pmatrix}
        \end{align}
        and
        \begin{align}
        \label{eq:G0def}
            \Gz = (\mathbf{\nabla_{\theta}Z})^\textrm{T} (\mathbf{\nabla_{\theta}Z}) & =      \begin{pmatrix}
\fraco{\partial\bm{Z}}{\partial \lambda} \cdot \fraco{\partial\bm{Z}}{\partial \lambda} & \fraco{\partial\bm{Z}}{\partial \lambda} \cdot \fraco{\partial\bm{Z}}{\partial \theta}\\[10 pt]
\fraco{\partial\bm{Z}}{\partial \theta} \cdot \fraco{\partial\bm{Z}}{\partial \lambda} & \fraco{\partial\bm{Z}}{\partial \theta} \cdot \fraco{\partial\bm{Z}}{\partial \theta}\\
        \end{pmatrix}.
  \end{align} 
The Green-Lagrange strain energy tensor is defined as
\begin{equation}
\label{eq:GLtensor}
\ten{\gamma} = \frac{1}{2}(\G-\Gz).
\end{equation}
In Cartesian coordinates, the Cauchy-Green deformation tensor is given by $\ten{C}=\ten{A}^\textrm{T}\ten{A}$, where $\ten{A}$ is the surface displacement gradient, $\ten{A}=\fraco{d\bm{X}}{d\bm{Z}}$ \cite{french}. As in \cite{fai} we consider a related tensor, 
\begin{equation}
\label{eq:Cten}
\ten{C}=\G \Gz^{-1}.
\end{equation}
In hyperelastic materials, the energy density functional $W$ is defined in terms of two deformation invariants, $I_1$ and $I_2$, where $I_1=$tr$(\C) - 2$ and $I_2=$det$(\C)-1$. As noted in \cite{french, deven}, the neo-Hookean energy density can be represented as a function of the invariants, i.e. $W(I_1,I_2)$, or as a function of the deformation gradient, $\mat{\nabla_{\theta}X}$. The total energy of the surface is given by an integral of the energy density
  \begin{align}
  \label{eq:energyintegral}
  \mathcal{E}(\vec{X})=\int_{\Gamma_0} W(I_1,I_2) \, dA_0 = \int_{\Gamma_0} W(\mat{\nabla_{\theta}X}) \, dA_0 = \int_{\Gamma_0} W(\mat{\nabla_{\theta}X}) \sqrt{\det{\Gz}} \, d\lambda \, d\theta,
  \end{align}
  where $\Gamma_0$ is the undeformed shell surface, and the area differential $dA_0$ denotes integration over the \textit{reference} configuration. The final equality comes from  the fact that $\sqrt{\det{\Gz}} = \norm{\vec{n_0}}$, where $\vec{n_0} = \fraco{\partial \vec{Z}}{\partial \lambda} \times \fraco{\partial \vec{Z}}{\partial \theta}$ is the normal (not necessarily a unit vector) to the surface in the reference configuration, as shown in \cite{fai}. We note that the factor $\sqrt{\det{\Gz}}$ is a function of the coordinates $\lambda$ and $\theta$, but \textit{not} of the current configuration $\vec{X}$ or its deformation gradient $\mat{\nabla_{\theta} X}$. 
  
In order to compute the Lagrangian force density Eq.\ \eqref{eq:forcefromvar},  consider a deformation by $\delta \vec{X}$. Evaluating Eq.\ \eqref{eq:vardef}, we have
\begin{align}
\phi_{\delta \bm{X}}(\alpha) = \mathcal{E}(\vec{X} + \alpha \delta \vec{X}) = \int_{\Gamma_0} W(\mat{\nabla_{\theta}}(\vec{X}+ \alpha \delta \vec{X})) \sqrt{\det{\Gz}} \, d\lambda \, d\theta.
\end{align}
To compute the variational derivative, we follow 
the derivation from \cite{deven} and first apply the chain rule for $\mat{\nabla_{\theta} X}$ (note that the area weight $\sqrt{\det{\Gz}}$ is unaffected),
\begin{equation}
\label{eq:chainruleeq}
\frac{d \phi_{\delta \bm{X}}(\alpha)}{d\alpha}\Bigg|_{\alpha=0} = \sum_{i=1}^3 \sum_{j=1}^2 {\int_{\Gamma_0} \frac{\partial W}{\partial (\mat{\nabla_{\theta} X})_{ij}} \frac{\partial \, \delta X_i}{\partial q_j}
\sqrt{\det{\Gz}} \, d\lambda \, d\theta}.
\end{equation}
The index $i$ describes the Cartesian coordinates (indexed from 1 to 3), and the index  $j$ describes the curvilinear coordinates $\lambda$ and $\theta$ (indexed from 1 to 2). We note that  
the first Piola-Kirchhoff stress tensor is 
\begin{equation}
\ten{P} = \fraco{\partial W}{\partial (\mat{\nabla_{\theta} X})}.
\end{equation}
Applying integration by parts on each of the coordinates $\lambda$ and $\theta$, we have
\begin{equation}
\frac{d \phi_{\delta \bm{X}}(\alpha)}{d\alpha}\Bigg|_{\alpha=0} = -\sum_{i=1}^3 \sum_{j=1}^2 {\int_{\Gamma_0} \frac{\partial (\ten{P}_{ij} \sqrt{\det{\Gz}})}{\partial q_j} \delta X_i
\, d\lambda \, d\theta}.
\end{equation}
The boundary terms evaluate to zero under the assumption that the surface is closed.
We obtain an area integral by multiplying and dividing by $\sqrt{\det{\Gz}}$,
\begin{equation}
\label{eq:frechsum}
\frac{d \phi_{\delta \bm{X}}(\alpha)}{d\alpha}\Bigg|_{\alpha=0} = -\sum_{i=1}^3 \sum_{j=1}^2 {\int_{\Gamma_0} \frac{1}{\sqrt{\det{\Gz}}}\frac{\partial (\ten{P}_{ij} \sqrt{\det{\Gz}})}{\partial q_j} \, \delta X_i \, 
dA_0}.
\end{equation}
Using the surface divergence operator in the undeformed state
\begin{equation}
\label{eq:surfacediv}
\nabla_{\theta}\,  \cdot \, = \frac{1}{\sqrt{\Det{\G_0}}} \left(\frac{\partial}{\partial \lambda} \sqrt{\Det{\G_0}} \textrm{ } \vec{e_1} + \frac{\partial}{\partial \theta} \sqrt{\Det{\G_0}} \textrm{ } \vec{e_2}\right) \cdot,
\end{equation}
we can write the variational derivative from Eq.\ \eqref{eq:forcefromvar} as
\begin{equation}
\label{eq:deldotP}
\frac{\delta{\mathcal{E}}}{\delta \bm{X}} = -\nabla_{\theta} \cdot \ten{P}^\textrm{T}.
\end{equation}
The above notation differs from \cite{deven} and instead follows \cite{hjelm}; it relies on a traditional tensor expansion in which $\ten{P}$ is written in terms of base dyads as $\displaystyle{\sum_{i,j} \ten{P}_{ij} \vec{e_i} \vec{e_j}}$. As shown in Appendix A, using $\nabla_{\theta}$ from Eq.\ \eqref{eq:surfacediv} and carrying out the dot product via tensor expansion in Eq.\ \eqref{eq:deldotP} yields the product in Eq.\ \eqref{eq:frechsum}. The force density per unit undeformed configuration in Eq.\ \eqref{eq:forcefromvar} is then  
\begin{equation}
\label{eq:feqdivdotPt}
\vec{F} = -\frac{\delta{\mathcal{E}}}{\delta \bm{X}} = \nabla_{\theta} \cdot \ten{P}^\textrm{T}.
\end{equation}
We point out the analogy with \cite{deven}, where the surface divergence operator in Cartesian coordinates has been replaced by the curvilinear surface divergence. Although we have obtained a representation for the force, it is in practice quite cumbersome to compute $\ten{P}=\fraco{\partial W}{\partial (\mat{\nabla_{\theta} X})}$ by directly differentiating the energy density (however, it is possible for some simple forms of $W$). 
We therefore expand the matrix $\ten{P}^{\textrm{T}}$ using the chain rule with the Green-Lagrange strain tensor, $\gamma$, as an intermediate:
  \begin{equation}
   \ten{P}^\textrm{T} = \frac{\partial{W}}{\partial{(\mathbf{\nabla_{\theta}X})^\textrm{T}}} = \frac{\partial W}{\partial \gamma} \cddot \frac{\partial \gamma}{\partial (\mathbf{\nabla_{\theta}X})^\textrm{T}},
   \end{equation}
where the double dot product $\cddot$ arises from the product of a fourth order tensor $\fraco{\partial \gamma}{\partial (\mathbf{\nabla_{\theta}X})^\textrm{T}}$ and second order tensor $\fraco{\partial W}{\partial \gamma}$  to form the second order tensor $\ten{P}^\textrm{T}$ (see Appendix A for a detailed expansion). We denote the symmetric second order tensor $\fraco{\partial W}{\partial \gamma} = \ten{S}$, where $\ten{S}$ is the second Piola-Kirchhoff stress tensor \cite{pozrikidis2001interfacial}. We show in Appendix A that
\begin{equation}
\label{eq:tenchain}
\ten{P}^{\textrm{T}}=\frac{\partial{W}}{\partial{(\mat{\nabla_{\theta}X})^\textrm{T}}} = \frac{\partial W}{\partial \gamma} \cddot \frac{\partial \gamma}{\partial (\mat{\nabla_{\theta}X})^\textrm{T}} = \ten{S} \cddot \frac{\partial \gamma}{\partial (\mat{\nabla_{\theta}X})^\textrm{T}} = \ten{S} \cdot (\mat{\nabla_{\theta}X})^\textrm{T}.
\end{equation}
By the symmetry of $\gamma$ (and $\ten{S}$), the transpose of both sides in Eq.\ \eqref{eq:tenchain} yields
\begin{equation}
\label{eq:PK12}
\ten{P}=(\mat{\nabla_{\theta}X}) \cdot \ten{S}^\textrm{T} = (\mat{\nabla_{\theta}X}) \cdot \ten{S}.
\end{equation}
In Eq.\ \eqref{eq:PK12}, we recover the correspondence between the first and second Piola-Kirchhoff stress tensors \cite{ogden1997non}.

To compute the tensor $\ten{S}$ in terms of the coordinates $\lambda$ and $\theta$, we begin by observing that in Eq.\ \eqref{eq:Cten}, $\C = \ten{I}_2 + 2\gamma \Gz^{-1}$, where $\ten{I}_2$ is the rank 2 identity tensor. Applying the chain rule and differentiating $\ten{C}$ directly with respect to $\gamma$ yields
\begin{equation}
\label{eq:sexpansion}
\ten{S} = \frac{\partial W}{\partial \gamma} = \frac{\partial W}{\partial \C} \cddot \frac{\partial \C}{\partial \gamma} = 2\frac{\partial W}{\partial \C} \cdot \Gz^{-1}.
\end{equation}
Computing $\fraco{\partial W}{\partial \C}$ by the chain rule, we have
\begin{equation}
\label{eq:dWdC}
\frac{\partial W}{\partial \C} = \frac{\partial W}{\partial I_1}\frac{\partial I_1}{\partial \C} + \frac{\partial W}{\partial I_2}\frac{\partial I_2}{\partial \C},
\end{equation}
where each $\fraco{\partial I_j}{\partial \C}$ is a second order tensor obtained from taking the derivative of each invariant $I_j$ with respect to each element of $\C$. We show in Appendix A that 
\begin{gather}
\label{eq:dI1dC}
\frac{\partial{I_1}}{\partial \C} = \ten{I}_2 \hspace{0.25cm}
\end{gather}
and
\begin{gather}
\label{eq:dI2dC}
\frac{\partial{I_2}}{\partial \C} = \left(\Det{\C}\right)\left(\G^{-1}\Gz\right).
\end{gather}
Substituting Eqs.\ \eqref{eq:dI1dC} and \eqref{eq:dI2dC} into Eqs. \eqref{eq:sexpansion} and \eqref{eq:dWdC},
\begin{gather}
\ten{S}  = 2\frac{\partial W}{\partial \C} \cdot \Gz^{-1} = 2\left(\frac{\partial W}{\partial I_1}\frac{\partial I_1}{\partial \C} + \frac{\partial W}{\partial I_2}\frac{\partial I_2}{\partial \C} \right) \cdot \Gz^{-1},\\[6 pt]
\ten{S} =  2\left(\frac{\partial W}{\partial I_1}\ten{I}_2 + \frac{\partial W}{\partial I_2}\left(\Det{\C}\right)\left(\G^{-1}\Gz\right) \right) \cdot \Gz^{-1},\\[6 pt] 
\ten{S} = 2\frac{\partial W}{\partial I_1}\Gz^{-1} + 2\frac{\partial W}{\partial I_2}\left(\Det{\C}\right) \G^{-1}.
\end{gather}
Since the constitutive law is often expressed in terms of the invariants $I_1$ and $I_2$, this formulation allows us to find $\vec{F}$ with respect to $\lambda$ and $\theta$ using Eq.\ \eqref{eq:PK12}.

Returning to the computation of force density as
$\vec{F}=\nabla_\theta \cdot \ten{P}^\textrm{T}$ from Eq.\ \eqref{eq:feqdivdotPt}, we apply the surface divergence operator in the undeformed state in Eq.\ \eqref{eq:surfacediv}. Then the force density can be expressed as
\begin{equation}
\label{eq:vdforce}
\vec{F} =  \frac{1}{\sqrt{\Det{\G_0}}} \left(\frac{\partial}{\partial \lambda} \sqrt{\Det{\G_0}} \left(\ten{S}_{11}\frac{\partial \vec{X}}{\partial \lambda} + \ten{S}_{12}\frac{\partial \vec{X}}{\partial \theta} \right) + \frac{\partial}{\partial \theta} \sqrt{\Det{\G_0}}\left(\ten{S}_{21}\frac{\partial \vec{X}}{\partial \lambda} + \ten{S}_{22}\frac{\partial \vec{X}}{\partial \theta} \right)\right).
\end{equation}
Further details are located in Appendix A. It is also possible to arrive at Eq.\ \eqref{eq:vdforce} by expanding the tensor $\mathcal{P}$ in Eq.\ \eqref{eq:chainruleeq} prior to integrating by parts.

%
%
%
%
%
\subsection{From force density to force}
%
%
%
The method developed in Section \ref{sec:vdderivation} allows us to find the force \textit{density} per unit reference configuration at any point $\vec{X}$. However, because the IB method requires the \textit{force} to be spread onto the fluid, it remains to compute the area weights associated with each point on the surface. While some applications do not require precise area weights, discretized structures in Stokes flow must have the correct area weights so that the discretized integral of force over the surface is zero. Because Eq.\ \eqref{eq:vdforce} is a translation invariant force density function, it can be shown \cite{peskinActa} that its integral is zero in the continuous sense. The discretization of the force integral must preserve this property for structures immersed in Stokes flow. In our formulation, this condition becomes one on the surface quadrature weights. If accurate surface quadrature weights can be computed, the force density will integrate to zero discretely (within some quadrature error); the reason is that the integral of the force generated by a
translation invariant hyperelastic energy over a closed surface is
known to be zero \cite{peskinActa}.

Let $n$ be the number of points at which the force is evaluated and consider the discretized integral of the force over the reference configuration
\begin{equation}
\label{eq:discintegral}
\int_{\Gamma_0} \vec{F} \, dA_0 \approx \sum_{i=1}^n \vec{F(\vec{X^i})}\, (\Delta A_0)_i.
\end{equation}
The factors $(\Delta A_0)_i$ are the appropriate area weights. Our goal is to determine the vector of area weights  $\vec{\Delta A_0} = \vec{\Delta A_0}(\vec{Z})$ (of length $n$) such that the integral in Eq.\ \eqref{eq:discintegral} is as close to exact as possible; that is, we are computing quadrature weights for integration over the undeformed configuration. Since the integral is over the reference configuration, the weights $\vec{\Delta A_0}$ only need to be computed once for the undeformed shape. In our applications, the undeformed shape is usually a sphere of radius 1, but the method we outline here can be used to compute the area weights in any reference configuration.

\subsubsection{Weights on a sphere of radius 1}
\label{sec:wtsspherer1}
Beginning with  the unit sphere $\mathbb{S}^2$, we use a method from \cite{graf} to compute nonnegative weights $\vec{\Delta A_0}$. Note that the nonnegativity constraint is necessary due to the physics of the problem; the force and force densities should point in the same directions, so the area weights should be nonnegative.

Let $n$ be the number of points at which the force function from Section \ref{sec:forcevariations} is evaluated (we refer to these as \textit{evaluation points}.) Following \cite{graf}, we compute nonnegative quadrature weights $(\Delta A_0)_i$ such that
\begin{equation}
\label{eq:discinequals}
\int \limits_{\mathbb{S}^2} f(\vec{X}) \, dA_0 = \sum_{i=1}^{n} f(\vec{X^i}) (\Delta A_0)_i 
\end{equation}
is exact for all functions, $f \in \Pi_N=\textrm{span}\left\{Y_{\ell}^k: \ell=0, 1, \dots N, k=-\ell, \dots,0, \dots \ell \right\}$. Here $\Pi_N$ is the $(N+1)^2$ dimensional space  of spherical harmonics up to degree $N$, defined as
\begin{equation}
\label{eq:Sharmonics}
Y_{\ell}^k(\lambda,\theta)=\begin{cases}
\sqrt{\fraco{(2\ell+1)}{4\pi}\fraco{(\ell-k)!}{(\ell+k)!}} P_{\ell}^k(\sin{\theta})\cos(k\lambda)\qquad & k \geq 0\\[8 pt]
\sqrt{\fraco{(2\ell+1)}{4\pi}\fraco{(\ell+k)!}{(\ell-k)!}} P_{\ell}^{-k}(\sin{\theta})\sin(-k\lambda)\qquad & k < 0,
\end{cases}
\end{equation}
where $P_{\ell}^k$ is an associated Legendre function of degree $\ell$ and order $k$. Let $f(\vec{X})=f(\lambda,\theta) = Y_0^0(\lambda,\theta) Y_\ell^k(\lambda,\theta)$. Then to satisfy Eq.\ \eqref{eq:discinequals},
\begin{align}
\int \limits_{\mathbb{S}^2} f(\vec{X}) \,  d\vec{x} & = \iint \limits_{\mathbb{S}^2} Y_0^0(\lambda,\theta)  Y_\ell^k(\lambda,\theta)\, dA_0\\[6 pt]
& = \sum_{i=1}^n Y_0^0(\lambda_i,\theta_i)  Y_\ell^k(\lambda_i,\theta_i) (\Delta A_0)_i\\
\label{eq:lasteq}
& = \frac{1}{\sqrt{4\pi}}\sum_{i=1}^n   Y_\ell^k(\lambda_i,\theta_i) (\Delta A_0)_i.
\end{align}
The last equality is due to the fact that $Y_0^0 = \fraco{1}{\sqrt{4\pi}}$. Since spherical harmonics form an orthonormal set, 
\begin{equation}
\label{eq:SHortho}
\int \limits_{\mathbb{S}^2} f(\vec{X})\, d\vec{x} = \iint \limits_{\mathbb{S}^2} Y_0^0(\lambda,\theta)  Y_\ell^k(\lambda,\theta) \, dA_0 = \begin{cases}
1 \quad \ell=0\\
0 \quad \ell \neq 0.
\end{cases}
\end{equation}
Equating Eqs.\ \eqref{eq:lasteq} and \eqref{eq:SHortho}, we have that 
\begin{align}
\label{eq:sysforwts}
\sum_{i=1}^n   Y_\ell^k(\lambda_i,\theta_i) (\Delta A_0)_i = \begin{cases}
\sqrt{4\pi} \quad & \ell=0\\
0 \quad & \ell \neq 0.
\end{cases}
\end{align}
Let the matrix
\begin{equation}
\label{eq:matY}
\mat{Y}=\begin{pmatrix}
Y_0^0(\lambda_1, \theta_1) & Y_1^0(\lambda_1, \theta_1) & Y_1^{-1}(\lambda_1, \theta_1) &  Y_1^{1}(\lambda_1, \theta_1) & \dots \\[4 pt]
Y_0^0(\lambda_2, \theta_2) & Y_1^0(\lambda_2, \theta_2) & Y_1^{-1}(\lambda_2, \theta_2) &  Y_1^{1}(\lambda_2, \theta_2) & \dots \\[4 pt]
\vdots & \vdots & \vdots & \vdots & \vdots\\[4 pt]
Y_0^0(\lambda_n, \theta_n) & Y_1^0(\lambda_n, \theta_n) & Y_1^{-1}(\lambda_n, \theta_n) &  Y_1^{1}(\lambda_n, \theta_n) & \dots
\end{pmatrix} \in \mathbb{C}^{n \times (N+1)^2}.
\end{equation}
Then in Eq.\ \eqref{eq:sysforwts} we seek the nonnegative least squares solution to
\begin{equation}
\label{eq:wtseq}
\mat{Y^*}\vec{\Delta A_0} = \sqrt{4\pi}\vec{e_1},
\end{equation}
where $\vec{e_1}=(1\quad 0 \quad \dots \quad 0)^{\textrm{T}}$ is an $(N+1)^2$ dimensional vector.
The nonnegative least squares solution for $\vec{\Delta A_0}$ gives the area weights for Eq.\ \eqref{eq:discinequals}. 

To solve Eq.\ \eqref{eq:wtseq}  for $\vec{\Delta{A}_0}$, we use the available implementation from \cite{nnls} to obtain the nonnegative area weights $\vec{\Delta {A}_0}$ for a unit sphere. Although any point set can be used in Eq.\ \eqref{eq:wtseq}, we let $n =(N+1)^2$ and choose evaluation points that maximize the determinant of the Gram matrix $\mathbf{YY^*}$. For these ``maximal determinant'' (MD) points, which are discussed in \cite{womersley2001good, sloan2004extremal} and can be downloaded from \cite{MDpts}, it has been conjectured that the quadrature weights satisfying the equality in Eq.\ \eqref{eq:wtseq} exist and are positive for all $N$ \cite{sloan2004extremal}.

Choosing weights that satisfy Eq.\ \eqref{eq:wtseq} ensures that the integral in Eq.\ \eqref{eq:discinequals} is exact for functions in $\Pi_N$, the space of spherical harmonics up to degree $N$. However, the calculated weights can still be used to integrate any function on the sphere. In fact, the force function in Eq.\ \eqref{eq:vdforce} is not necessarily in $\Pi_N$ since it involves derivatives of arbitrary position functions. A quadrature error is therefore introduced. In \cite{sloan2004extremal}, it is shown that the quadrature error, $\epsilon$, in a certain Hilbert space is bounded by
\begin{equation}
\epsilon \leq \fraco{\sqrt{4\pi}}{(N+1)^{1/2}}=\fraco{\sqrt{4\pi}}{n^{1/4}},
\end{equation}
where the last equality holds only in the case that $n=(N+1)^2$. For MD points, the measured worst case error is much better; it is $\mathcal{O}(n^{-3/4})$ (see \cite{sloan2004extremal}, Table 2). Our observations in Section \ref{sec:forceaccuracy} confirm this, as the actual quadrature error in integrating the force function was found to be several orders of magnitude lower than the theoretical bound, even for a highly irregular shape.

\subsubsection{Weights on an arbitrary reference configuration}
\textbf{}We use the formulation from \cite{grady} to determine the quadrature weights on an arbitrary configuration. Let $\mathbb{L}$ denote a space diffeomorphic to $\mathbb{S}^2$, and let $\displaystyle{\vec{Z}_{\mathbb{S}^2}(\lambda,\theta)}$ and $\displaystyle{\vec{Z}_{\mathbb{L}}(\lambda,\theta)}$ be the mappings from Eq.\ \eqref{eq:surfacereps} that describe $(\hat{X},\hat{Y},\hat{Z})$ from $(\lambda,\theta)$ on each space. Then, as shown in \cite{grady}, the surface elements for each configuration are given by
\begin{align}
\label{eq:volcurr}
dA_{\mathbb{L}}& = \sqrt{\det{\left(\Gz\left(\vec{Z}_\mathbb{L}\right)\right)}} \, d\lambda \, d\theta
\end{align}
and
\begin{align}
\label{eq:volr1}
dA_{\mathbb{S}^2}& = \sqrt{\det{\left(\Gz\left(\vec{Z}_{\mathbb{S}^2}\right)\right)}} \, d\lambda \, d\theta,
\end{align}
where $\Gz$ is evaluated for each surface according to Eq.\ \eqref{eq:G0def}. From Eqs.\ \eqref{eq:volcurr} and \eqref{eq:volr1}, it is clear that
\begin{equation}
dA_{\mathbb{L}}= \sqrt{\frac{\det{\left(\Gz\left(\vec{Z}_\mathbb{L}\right)\right)}}{\det{\left(\Gz\left(\vec{Z}_{\mathbb{S}^2}\right)\right)}}} \, dA_{\mathbb{S}^2}.
\end{equation}
Thus the analogous integral to Eq.\ \eqref{eq:discintegral} is
\begin{equation}
\label{eq:arbitintegral}
\int_{\mathbb{L}} f(\vec{X}) \, {dA}_\mathbb{L} = \int_{\mathbb{S}^2} f(\vec{X}) \, \sqrt{\frac{\det{\left(\Gz\left(\vec{Z}_\mathbb{L}\right)\right)}}{\det{\left(\Gz\left(\vec{Z}_{\mathbb{S}^2}\right)\right)}}}\, dA_{\mathbb{S}^2} = \sum_{i=1}^{n} f(\vec{X^i}) \sqrt{\frac{\det{\left(\Gz\left(\vec{Z}_\mathbb{L}\right)\right)}}{\det{\left(\Gz\left(\vec{Z}_{\mathbb{S}^2}\right)\right)}}} \, \left(\Delta {A}_0\right)_i.
\end{equation}
According to Eq.\ \eqref{eq:arbitintegral}, the weights for the arbitrary reference configuration are given by
\begin{equation}
\label{eq:arbitwts}
\vec{\Delta {A}_\mathbb{L}} = \sqrt{\frac{\det{\left(\Gz\left(\vec{Z}_\mathbb{L}\right)\right)}}{\det{\left(\Gz\left(\vec{Z}_{\mathbb{S}^2}\right)\right)}}} \vec{\Delta {A}_{\mathbb{S}^2}}.
\end{equation}
For example, if the reference configuration is a sphere of radius $r$, then $\det{\left(\Gz\left(\vec{Z}_\mathbb{L}\right)\right)} = r^4\det{\left(\Gz\left(\vec{Z}_{\mathbb{S}^2}\right)\right)}$, so $\vec{\Delta {A}_\mathbb{L}} = r^2 \vec{\Delta {A}_{\mathbb{S}^2}}$, as expected since the surface area of a sphere scales by $r^2$. 

With the area weights from Eq.\ \eqref{eq:arbitwts}, we have finally found the \textit{force}, $\vec{\widehat{F}}$ at point $i$ from the force density \vec{F} as a function of position in $(\lambda,\theta)$,
\begin{equation}
\label{eq:forcefromforced}
\vec{\widehat{F}}_i(\lambda,\theta) = \vec{F}_i(\lambda,\theta) (\Delta A_0)_i.
\end{equation}

\subsection{The linear discrete surface method (LDSM)}
\label{sec:LDSM}
In LDSM (first introduced in \cite{kraus1996fluid}), the surface of the shell is triangulated and the integral for energy is discretized as such. Consider an energy functional from Eq.\ \eqref{eq:energyfunc} of the form
\begin{equation}
  \mathcal{E}=\int_{\Gamma_0} W(I_1, I_2) \, dA_0.
\end{equation}
The above integral can be discretized into a sum over the surface triangles as follows,
\begin{equation}
\mathcal{E}=\sum_{\textrm{tri}} \int_{\textrm{tri}} W(I_1, I_2) \, dA_0.
\end{equation}
LDSM uses barycentric coordinates ($\gamma_1, \gamma_2, \gamma_3$) to represent a deformed triangle $t$ with vertices $\vec{X}^1, \vec{X}^2, \vec{X}^3$ in the deformed configuration that corresponds to a reference triangle $s$ with vertices $\vec{Z}^1, \vec{Z}^2, \vec{Z}^3$ in the undeformed configuration. Any point within the triangle can be written as
\begin{align}
\vec{X} & =\gamma_1 \vec{X}^1 + \gamma_2 \vec{X}^2 + \gamma_3 \vec{X}^3 
\end{align}
and 
\begin{align}
\vec{Z} & =\gamma_1 \vec{Z}^1 + \gamma_2 \vec{Z}^2 + \gamma_3 \vec{Z}^3,
\end{align}
where $\gamma_1+\gamma_2+\gamma_3=1$ and $\gamma_1,\gamma_2,\gamma_3 \geq 0$. Using the relation $\gamma_3= 1 - \gamma_1-\gamma_2$, the surface can be represented in terms of two material coordinates (note that, for $\gamma_3 \geq 0$ to be satisfied, $\gamma_1 \leq 1-\gamma_2$). As in the continuous case, we define the matrices
  \begin{align}
     \mat{\nabla_{q}X}& =\begin{pmatrix}
        \fraco{\partial\vec{X}}{\partial \gamma_1} & \fraco{\partial\vec{X}}{\partial \gamma_2}
        \end{pmatrix} = \begin{pmatrix}
        \vec{X}^1 - \vec{X}^3 & \vec{X}^2 - \vec{X}^3
        \end{pmatrix}  
        \end{align}
and
        \begin{align}
        \mat{\nabla_{q}Z}& =\begin{pmatrix}
        \fraco{\partial\vec{Z}}{\partial \gamma_1} & \fraco{\partial\vec{Z}}{\partial \gamma_2}
        \end{pmatrix} = \begin{pmatrix}
        \vec{Z}^1 - \vec{Z}^3 & \vec{Z}^2 - \vec{Z}^3
        \end{pmatrix}.
  \end{align}
Denote $\vec{X}^i - \vec{X}^j = \vec{X}^{ij}$ (and likewise for \vec{Z}). Then we can again define the metric tensors
\begin{align}
     \G = (\mathbf{\nabla_{q}X})^\textrm{T} (\mathbf{\nabla_{q}X}) & =      \begin{pmatrix}
\vec{X}^{13} \cdot \vec{X}^{13}  & \vec{X}^{13} \cdot \vec{X}^{23}\\[6 pt]
\vec{X}^{23} \cdot \vec{X}^{13}  & \vec{X}^{23} \cdot \vec{X}^{23}\\
        \end{pmatrix} 
        \end{align}
        and
        \begin{align}
            \Gz = (\mathbf{\nabla_{q}Z})^\textrm{T} (\mathbf{\nabla_{q}Z}) & =      \begin{pmatrix}
\vec{Z}^{13} \cdot \vec{Z}^{13}  & \vec{Z}^{13} \cdot \vec{Z}^{23}\\[6 pt]
\vec{Z}^{23} \cdot \vec{Z}^{13}  & \vec{Z}^{23} \cdot \vec{Z}^{23}\\
        \end{pmatrix}.
  \end{align} 

The tensors $\G$ and $\Gz$ are \textit{constant} on any triangle $t$ or $s$. Therefore, the tensor $\C = \G \Gz^{-1}$ and its invariants $I_1$ and $I_2$ are constant on $t$ and $s$. Thus any energy density functional defined in terms of the invariants is also constant on $t$ and $s$ and can be defined in terms of the vertices $\vec{X}^1, \vec{X}^2, \vec{X}^3,$ $\vec{Z}^1, \vec{Z}^2, \textrm{ and } \vec{Z}^3$. Denote the energy density of a triangle as $W(t,s)$ (it is a function of the deformed triangle $t$ and undeformed triangle $s$). The discretized energy is then
\begin{equation}
\mathcal{E} =\sum_{\textrm{tri}} \int_{s} W(I_1, I_2) \, dA_0 = \sum_{\textrm{tri}} \int_{s} W(t,s) (\Det{\Gz})^{1/2} \, d\gamma_1 \, d\gamma_2.
\end{equation}
Note that the integration is carried out over the \textit{reference} triangle $s$ and that $dA_0 = (\Det{\Gz})^{1/2} \, d\gamma_1 \, d\gamma_2$ since $(\Det{\Gz})^{1/2}$ is the magnitude of the normal to the reference triangle \cite{fai}. Evaluating the integral, we have
\begin{align}
\mathcal{E} & = \sum_{\textrm{tri}} W(t,s) (\Det{\Gz})^{1/2} \int_{s}  d\gamma_2 \, d\gamma_1, \\[6 pt]
\mathcal{E} & = \sum_{\textrm{tri}} W(t,s) (\Det{\Gz})^{1/2} \int_{0}^{1} \int_{0}^{1-\gamma_1}  d\gamma_2 \, d\gamma_1,\\[6 pt]
\mathcal{E} & = \frac{1}{2}\sum_{\textrm{tri}} W(t,s) (\Det{\Gz})^{1/2}.
\end{align}
From the above equation, note that $\Delta A_0(s) = \fraco{1}{2}(\Det{\Gz})^{1/2}$, which is expected \textit{a priori} since it is simply the formula for the area of the reference triangle $s$. Since the energy functional has already been integrated, the area weights at each point have already been accounted for, and the \textit{force} (not force density), $\vec{\widehat{F}}$, at a point $\vec{X}^i$ due to triangle $t$ can be computed by taking the derivative of the energy with respect to that point,
\begin{equation}
\vec{\widehat{F}}_i(t) = -\frac{\partial \mathcal{E}}{\partial \vec{X}^i}\left(t\right) = -\frac{\partial }{\partial \vec{X}^i}\left(\frac{1}{2} W(t,s) (\Det{\Gz})^{1/2}\right).
\end{equation}
The force at each point is then the sum of the force contributions from each triangle that includes the point, i.e.
\begin{equation}
\label{eq:faiforce}
\vec{\widehat{F}}_i\left(\vec{X}^i\right) = \sum_{t \ni \vec{X}^i}-\frac{\partial \mathcal{E}}{\partial \vec{X}^i}\left(t\right) = \sum_{t \ni \vec{X}^i} -\frac{\partial }{\partial \vec{X}^i}\left(\frac{1}{2} W(t,s) (\Det{\Gz})^{1/2}\right).
\end{equation}
It can be shown using vector calculus that the sum of the forces over the points in LDSM is zero. The force in Eq.\ \eqref{eq:faiforce} therefore satisfies the discrete integral constraint for Stokes flow. The force density per unit undeformed area is computed by dividing the force at each point by its associated area weight in the undeformed configuration. The area weight at point $i$ in LDSM is given by $1/3$ of the sum of the reference triangle areas that include point $i$. That is, 
\begin{equation}
\label{eq:faiwts}
\Delta (A_{L_0})_i = \frac{1}{3} \sum_{s \ni \vec{Z}^i} { ({\det{\Gz(s)}})^{1/2} }.
\end{equation}
Further details are located in \cite{fai}.
%
%
%
%
\subsection{Energy functionals}
\label{sec:energyfunctionals}
%
%
%
In this section, we describe energy functionals for neo-Hookean energy, surface tension, and bending (or  curvature) energy.
\subsubsection{Neo-Hookean energy}
We follow Evans \& Skalak's formulation for the neo-Hookean energy density \cite{ES} (for other possible neo-Hookean constitutive laws, see \cite{french}). In terms of the invariants,
\begin{align}
\label{eq:nhforce}
W_{NH}  & = G_s\left(\frac{I_1 + 2}{2\sqrt{I_2+1}} - 1 + \frac{A}{2}\left(I_2 + 2 - 2\sqrt{I_2+1}\right)\right)\\[6 pt]
& = G_s\left(\frac{\Tr{\C}}{2(\Det{\C})^{1/2}} - 1 + \frac{A}{2}\left(\Det{\C} + 1 - 2(\Det{\C})^{1/2}\right)\right)\\[6 pt]
& = G_s\left(\frac{\Tr{\C}}{2(\Det{\C})^{1/2}} - 1 + \frac{A}{2}\left((\Det{\C})^{1/2} -1 \right)^2\right)\\[6 pt]
\label{eq:faiform}
& = \frac{K}{2}\left((\Det{\C})^{1/2} -1 \right)^2 + \frac{G_s}{2}\left(\frac{\Tr{\C}}{(\Det{\C})^{1/2}} - 2\right) ,
\end{align}
where $K=AG_s$ is the area dilation modulus and $G_s$ is the shear modulus. Eq.\ \eqref{eq:faiform} gives the form used in \cite{fai}. We take derivatives of Eq.\ \eqref{eq:nhforce} with respect to $I_1$ and $I_2$ to find the tensors $\ten{P}$, $\ten{S}$, and force density $\vec{F}$ in Eq.\ \eqref{eq:vdforce}.

\subsubsection{Surface energy}
\label{sec:surfacetens}
The most common representation of ``surface energy'' is $\mathcal{E}_{ST} = \sigma \int dA$, where $dA$ is the area element in the \textit{current} configuration and $\sigma$ is the surface tension (with units of force/length) of the surface. For both the variational force density method and LDSM, $dA = (\Det{\G})^{1/2} \, d\lambda \, d\theta$, but this must be expressed with respect to the integral over the reference area to satisfy our method. Therefore we use Eq.\ \eqref{eq:arbitintegral} and \cite{grady} to write
\begin{equation}
\mathcal{E}_{ST} = \int_{\Gamma} \sigma \, dA = \int_{\Gamma_0} \sigma \left(\frac{\Det{\G}}{\Det{\Gz}}\right)^{1/2} \, dA_0 = \int_{\Gamma_0} \sigma \left(\Det{\C}\right)^{1/2} \, dA_0.
\end{equation}
We can then express the surface energy in the form of Eq.\ \eqref{eq:energyintegral}  with $W_{ST}=\sigma \left(\Det{\C}\right)^{1/2} = \sigma(I_2 + 1)^{1/2}$. 
Derivatives of this density with respect to $I_1$ and $I_2$ can easily be calculated to evaluate the first and second Piola-Kirchhoff stress tensors $\ten{P}$, $\ten{S}$ and force density $\vec{F}$ in Eq.\ \eqref{eq:vdforce}.

The variational derivative can be carried out over the current configuration as well to yield force density per unit deformed area \cite{shankar,capovilla2017elastic,pozrikidis2003modeling}. When this formulation, which is related directly to the mean curvature of the surface, is used, a transformation must be made to express the force per unit reference area.
\begin{equation}
\vec{F}_{cur}^{ST} \, dA= \sigma (2H)\, \hat{\vec{n}} \, dA =\sigma (2H)\, \hat{\vec{n}} \,  \left(\frac{\Det{\G}}{\Det{\Gz}}\right)^{1/2} \, dA_0.
\end{equation}
The force per unit undeformed configuration then is, 
\begin{equation}
\label{eq:stforce}
\vec{F}^{ST} = \sigma (2H)\, \hat{\vec{n}} \,  \left(\frac{\Det{\G}}{\Det{\Gz}}\right)^{1/2},
\end{equation}
where $\hat{\vec{n}}$ is the unit normal vector in the current configuration,
\begin{equation}
\label{eq:unitnorm}
\hat{\vec{n}}=\fraco{\fraco{\partial\bm{X}}{\partial \lambda} \times \fraco{\partial\bm{X}}{\partial \theta}}{\norm{\fraco{\partial\bm{X}}{\partial \lambda} \times \fraco{\partial\bm{X}}{\partial \theta}}}.
\end{equation}
$H$ is the mean curvature (one half the total curvature, or the average of the two principle curvatures) given by the formula
\begin{equation}
\label{eq:meancurve}
H  = \frac{eG - 2fF + gE}{2(EG - F^2)},
\end{equation}
where $E$, $F$, and $G$ are coefficients of the first fundamental form
\begin{equation}
E = \frac{\partial\vec{X}}{\partial \lambda} \cdot \frac{\partial\vec{X}}{\partial \lambda} \qquad
F = \frac{\partial\vec{X}}{\partial \lambda} \cdot \frac{\partial\vec{X}}{\partial \theta} \qquad
G = \frac{\partial\vec{X}}{\partial \theta} \cdot \frac{\partial\vec{X}}{\partial \theta},
\end{equation}
and $e$, $f$, and $g$ are coefficients of the second fundamental form,
\begin{equation}
\label{eq:secffs}
e = \frac{\partial^2 \vec{X}}{\partial \lambda^2} \cdot \hat{\vec{n}} \qquad
f = \frac{\partial^2 \vec{X}}{\partial \lambda \textrm{ } \partial \theta} \cdot\hat{\vec{n}} \qquad
g = \frac{\partial^2 \vec{X}}{\partial \theta^2} \cdot \hat{\vec{n}}.
\end{equation}
It can be shown using elementary vector calculus that Eq.\ \eqref{eq:stforce} is equivalent to Eq.\ \eqref{eq:vdforce} when $W_{ST}=\sigma(I_2 + 1)^{1/2}$ is used. We emphasize that for this simple case, we were able to take a surface energy given over the current configuration and transform it so we could integrate over the reference configuration and apply our method. 

More complex energy functionals given over the current configuration are not as easily transformed. In these cases, differential geometry can be used to find the variational derivatives and force per unit deformed area. These expressions can then be transformed to the undeformed configuration by multiplying by the ratio of area weights as in Eq.\ \eqref{eq:stforce}. Bending energy is an example of such an energy functional and is described in the next section.


\subsubsection{Bending energy}
\label{sec:brforce}
In this section, we describe the variational formulation of forces due to bending. It is worth noting that piecewise linear representations of the surface (such as LDSM) must use a method different from that of Section \ref{sec:LDSM}, where the energy is evaluated on each triangle, to compute forces due to bending. This is because bending energies are directly related to curvature, and in a piecewise linear surface representation the triangles have zero curvature. In the case of \cite{fai}, the method used to compute forces due to bending does not ensure that the bending force integrates to zero (despite the bending energy being translation invariant), and therefore cannot be used in the context of Stokes flow. Our method allows for the true bending energy functional to be used because we can properly resolve the curvature of the surface. In addition, the use of the correct quadrature weights ensures that the bending force integrates to zero over the fluid domain. 

Interestingly, variational principles have been used to derive forces due to bending rigidity for quite some time, with one early example from \cite{zhong1989bending}. These previous studies addressed the problem of determining the equilibrium shape of membranes experiencing bending forces by setting the force density equal to zero. The derivation of the force density due to bending begins with the energy functional
\begin{equation}
\label{eq:bendenergy}
\mathcal{E}_{bend}=k_{bend} \int_\Gamma (2H)^2 \, dA,
\end{equation}
where $k_{bend}$ is the bending rigidity and the energy is an integral of the total curvature ($2H$) over the \textit{deformed} configuration. The variational approach calls for the energy functional to be written in the form
\begin{equation}
\delta \mathcal{E}_{bend} =  \int_\Gamma -\vec{F}^{BR}_{cur} \cdot \delta \vec{X} \, dA.
\end{equation}
As shown in \cite{capovilla2017elastic, capovilla2004second}, arguments from differential geometry can be used to write the force density per unit deformed configuration as
\begin{equation}
\label{eq:brforcecur}
\vec{F}^{BR}_{cur}=-k_{bend}\left(4{\nabla^2_s} H + 8H^3 -4HR\right)\hat{\vec{n}},
\end{equation}
where $H$ is the mean curvature defined according to Eq.\ \eqref{eq:meancurve}, $\hat{\vec{n}}$ is the unit normal defined in Eq.\ \eqref{eq:unitnorm}, and ${\nabla^2_s}$ is the surface Laplacian in the deformed state, defined as
\begin{equation}
\label{eq:surfacelapcur}
{\nabla^2_s} = \sum_{i=1}^2 { \sum_{j=1}^2 {\fraco{1}{\sqrt{\det{\G}}} \fraco{\partial}{\partial q_i}\left(\sqrt{\det{\G}}\, \G^{-1}_{ij} \fraco{\partial}{\partial q_j}\right) }}.
\end{equation}
Finally, $R$ is the scalar curvature, or twice the Gaussian curvature, defined as
\begin{equation}
R = 2 \det{\Bigg(\G^{-1}\begin{pmatrix} e & f \\[4 pt]
f & g \end{pmatrix}\Bigg)},
\end{equation}
where the fundamental forms $e$, $f$, and $g$ are defined in Eq.\ \eqref{eq:secffs}. 
In a similar manner to that of Eq.\ \eqref{eq:stforce}, Eq.\ \eqref{eq:brforcecur} can be written per unit undeformed area as
\begin{equation}
\label{fbrundef}
\vec{F}^{BR}=-k_{bend}\left(\frac{\Det{\G}}{\Det{\Gz}}\right)^{1/2}\left(4{\nabla^2_s} H + 8H^3 -4HR\right)\hat{\vec{n}} ,
\end{equation}
which gives the force due to bending per unit undeformed area. 
%
%
\section{Numerical Formulation}
\label{sec:numform}
%
%
In this section, we describe the discretization of the surface and fluid grids. We then describe the fluid solver and temporal update for dynamic simulations.
\subsection{Lagrangian surface and force discretization}
For the variational force density method, some mapping must be chosen for $\vec{X}$ in Eq.\ \eqref{eq:surfacereps}. Our method allows any mapping that defines a closed surface to be used, including interpolants composed of Lagrange, spherical harmonic, or radial basis functions. We choose to use a spherical harmonic interpolant for the maps $X$, $Y$, and $Z$ in Eq.\ \eqref{eq:surfacereps} (the reference configuration $\vec{Z}$ is generally taken to be the unit sphere unless otherwise specified). Each interpolant is computed in the same manner and is also described in \cite{shankar}. Suppose we have a set of $m = (M+1)^2$ points in $(x,y,z)$ space. Then we use a spherical harmonic interpolant of order $M$ to approximate the function $\vec{X}(\lambda,\theta)$ in Eq.\ \eqref{eq:surfacereps}. For example,
\begin{equation}
\label{eq:shintetp}
X(\lambda,\theta) = \sum_{\ell=0}^{M}\left(\sum_{k=-\ell}^{-1} c^x_{\ell,k} Y_\ell^k(\lambda,\theta) + \sum_{k=0}^{\ell} c^x_{\ell,k} Y_\ell^k(\lambda,\theta) \right),
\end{equation}
where the spherical harmonics are as defined in Eq.\ \eqref{eq:Sharmonics} and the other interpolants $Y$, $Z$, $\hat{X}$, $\hat{Y}$, and $\hat{Z}$ are defined similarly. In \cite{shankar}, it is shown that Eq.\ \eqref{eq:shintetp} corresponds to a linear system
\begin{equation}
\label{eq:systeminterp}
\mathbf{Y}\vec{c^x} = 
\begin{pmatrix}
Y_0^0(\lambda_1, \theta_1) & Y_1^0(\lambda_1, \theta_1) & Y_1^{-1}(\lambda_1, \theta_1) &  Y_1^{1}(\lambda_1, \theta_1) & \dots \\[4 pt]
Y_0^0(\lambda_2, \theta_2) & Y_1^0(\lambda_2, \theta_2) & Y_1^{-1}(\lambda_2, \theta_2) &  Y_1^{1}(\lambda_2, \theta_2) & \dots \\[4 pt]
\vdots & \vdots & \vdots & \vdots & \vdots\\[4 pt]
Y_0^0(\lambda_m, \theta_m) & Y_1^0(\lambda_m, \theta_m) & Y_1^{-1}(\lambda_m, \theta_m) &  Y_1^{1}(\lambda_m, \theta_m) & \dots
\end{pmatrix}
\begin{pmatrix}
c_1^x\\[4 pt]
c_2^x\\[4 pt]
\vdots\\[4 pt]
c_m^x
\end{pmatrix}
=\begin{pmatrix}
X_1\\[4 pt]
X_2\\[4 pt]
\vdots\\[4 pt]
X_m
\end{pmatrix},
\end{equation}
where $m = (M+1)^2$ exactly and $X_1, X_2, \dots, X_m$ are the $x$ coordinates of the first, second, and $m$th points, respectively. Since the coordinates of the interpolation points are constant throughout the algorithm, the $LU$ factorization of the matrix $\mathbf{Y}$ in Eq.\ \eqref{eq:systeminterp} can be precomputed at the start of the algorithm. Then at each timestep, $\mathcal{O}(m^2)$ flops are required to solve Eq.\ \eqref{eq:systeminterp}.

For the variational force density method, we follow \cite{shankar} and define the $m$ points in $(x,y,z)$ space that are used to construct the interpolant as \textit{interpolation} points. Note that there are typically much fewer interpolation points than evaluation points ($ m \ll n$). Because the matrix in Eq.\ \eqref{eq:systeminterp} must be invertible, we use the maximal determinant points mentioned previously (downloaded from \cite{MDpts}) that maximize the determinant of the Gram matrix $\mathbf{YY^*}$. See \cite{womersley2001good} for more details on choosing points on a sphere for interpolation.

The interpolant from Eq.\ \eqref{eq:shintetp} is the map $\vec{X}(\lambda,\theta)$ in Eq.\ \eqref{eq:surfacereps}, and given a reference configuration $\vec{Z}(\lambda,\theta)$ we proceed in computing force densities from the variational method in Eq.\ \eqref{eq:vdforce} at a different set of $n$ maximal determinant evaluation points. We find the spherical coordinates of the evaluation points and use the function $\vec{X}$ to map them to our deformed configuration, where we evaluate the forces using Eq.\ \eqref{eq:vdforce} with the area weights from Eq.\ \eqref{eq:wtseq}. The calculation of the forces in  Eq.\ \eqref{eq:vdforce} requires the derivatives of $\vec{X}$ at each of the evaluation points. In practice, this is accomplished by multiplying matrices of spherical harmonic derivatives (e.g. $\fraco{\partial{\mathbf{Y}}}{\partial \lambda}$, where the derivatives are taken on each element) by the vector of function weights ($\vec{c^j}$). Each derivative is thus a matrix-vector multiplication problem and requires $\mathcal{O}(mn)$ operations, although this could easily be parallelized to reduce computational time. Once the derivatives are computed, $\mathcal{O}(n)$ operations are required to compute the forces in Eq.\ \eqref{eq:vdforce} from the derivative values (this step could also easily be done in parallel). The leading order term in the flop count for the force calculation in our method is therefore $\mathcal{O}(mn)$.

In LDSM, the force is computed directly at each of the evaluation points using Eq.\ \eqref{eq:faiforce} with the triangles from the surface discretization. Since the force is computed triangle by triangle and there are approximately $2n$ triangles, the entire force calculation for LDSM is $\mathcal{O}(n)$. 
We perform timing tests in Section \ref{sec:dynamics} to compare the computational costs of our method with LDSM in practice.

Each of the weights $(\Delta A_0)_i$ in Eq.\ \eqref{eq:wtseq} is computed for the $n$ evaluation points and multiplied by its respective force density to get a total force at each of the evaluation points. Since the continuous force density representation is exact for a given position function $\vec{X}$, the only errors introduced are in representing the surface by a continuously differentiable spherical harmonic interpolant and in using a quadrature rule to integrate the continuous force density function. A detailed study of the errors associated with a spherical harmonic representation can be found in \cite{shankar}.  An analytic estimate for the quadrature error is located in \cite{sloan2004extremal}. We also analyze the error in force computations and compare it to the error in LDSM in Section \ref{sec:forceaccuracy}.

\subsection{Eulerian discretization}
The fluid domain is discretized on a periodic three dimensional domain $[-L,L] \times [-L,L] \times [-L,L]$ with mesh spacing $\Delta x= \Delta y= \Delta z = 2L/\eta$, where $\eta$ is the grid size. We use a spectral fluid solver to compute the pressure and velocity at each time step in $\mathcal{O}(\eta^3 \log{\eta^3})$ operations \cite{PeyretBook}. The evaluation points are chosen so that approximately two Lagrangian points lie within an Eulerian cube of dimension $\Delta x^3$.
\subsection{Time update}
\label{sec:timeupdate}
We outline the time stepping procedure for the variational derivative method first. 
Given a set of $m$ interpolation points and $n$ evaluation points on the reference configuration, 
and the map $\vec{Z}(\lambda,\theta)$ in Eq.\ \eqref{eq:surfacereps}, the area weights $\vec{\Delta {A}_0}$ in Eq.\ \eqref{eq:wtseq}, and matrix $\mathbf{Y}$ in Eq.\ \eqref{eq:systeminterp} (and its $LU$ factorization) are precomputed. Then at the $k$-th time step,
\begin{enumerate}
\item The $m$ interpolation points in the current configuration $(x_i,y_i,z_i)$ are used to solve Eq.\ \eqref{eq:systeminterp} for the map 
$\vec{X}$ ($\mathcal{O}(m^2)$ operations).
\item The computed maps $\vec{X}$ and $\vec{Z}$ are used to find the force densities in Eq.\ \eqref{eq:vdforce} at each of 
the $n$ \textit{evaluation} points ($\mathcal{O}(mn)$ operations).
\item The forces are spread onto the Eulerian grid using the discrete delta function from \cite{peskinActa} 
with the area weights from Eq.\ \eqref{eq:wtseq} that ensure the spread force integrates to zero ($\mathcal{O}(4^3n)$ operations because we are using a 4 point discrete delta function).
\item The Stokes equations are solved on the Eulerian grid ($\mathcal{O}(\eta^3 \log{\eta^3})$ operations).
\item  The fluid velocity is interpolated back to the structure at the $m$ \textit{interpolation} points ($\mathcal{O}(4^3m)$ operations).
\item The interpolation points are updated with the local fluid velocity,
\begin{equation}
\vec{X}_{i}^{k+1} = \vec{X}^{k}_i + \Delta t \, \vec{U}^{k+1},
\end{equation}
where $\vec{U}^{k+1}$ indicates the interpolated fluid velocity from step 5 ($\mathcal{O}(m)$ operations).
\end{enumerate}
In LDSM, we work with only the set of $n$ evaluation points.
We begin with a triangulated surface and use the locations of the points in the current and reference configurations to compute the forces, then proceed with steps 3-6 as above, with the exception being that the structure velocity is computed at the same $n$ points where the force was evaluated. The interpolation of the velocity in LDSM therefore requires $\mathcal{O}(4^3n)$ operations instead of $\mathcal{O}(4^3m)$, and the forward Euler update is $\mathcal{O}(n)$ in LDSM instead of $\mathcal{O}(m)$ in our method. If $m \ll n$, our method results in fewer operations in steps 5 and 6.
%
%
%
%
%
%
\section{Results}
\label{sec:results}
In this section we test the accuracy of LDSM and the spherical harmonic/variational derivative (SHVD) approaches in the context of both force computations and dynamic immersed boundary simulations. We first consider the accuracy of the force computation on an ellipsoid and perturbed ellipsoid and proceed with dynamic simulations on each surface.
\subsection{Accuracy of force computations}
%
%
\label{sec:forceaccuracy}
We consider the accuracy of surface tension and neo-Hookean forces on two surfaces given by
\begin{gather}
\label{eq:obj1}
\textrm{Ellipsoid: } \vec{X} = \begin{pmatrix}
a \cos{\lambda} \cos{\theta}\\[4 pt]
b \sin{\lambda} \cos{\theta}\\[4 pt]
c \sin{\theta}
\end{pmatrix}
\end{gather}
and 
\begin{gather}
\label{eq:obj2}
\textrm{``Perturbed ellipsoid:'' } \vec{X} = V_f \begin{pmatrix}
a \left(1+\fraco{B}{5}\exp(-\sin{\theta})\right) \cos{\lambda}\cos{\theta}\\[8 pt]
b \left(1+B\exp(-\sin{\theta})\right) \sin{\lambda} \cos{\theta}\\[8 pt]
c \left(1+B\exp(-\sin{\theta})\right) \sin{\theta}
\end{pmatrix}.
\end{gather}
\begin{figure}
\centering     
\subfigure[Ellipsoid]{\label{fig:shapea}\includegraphics[width=75mm]{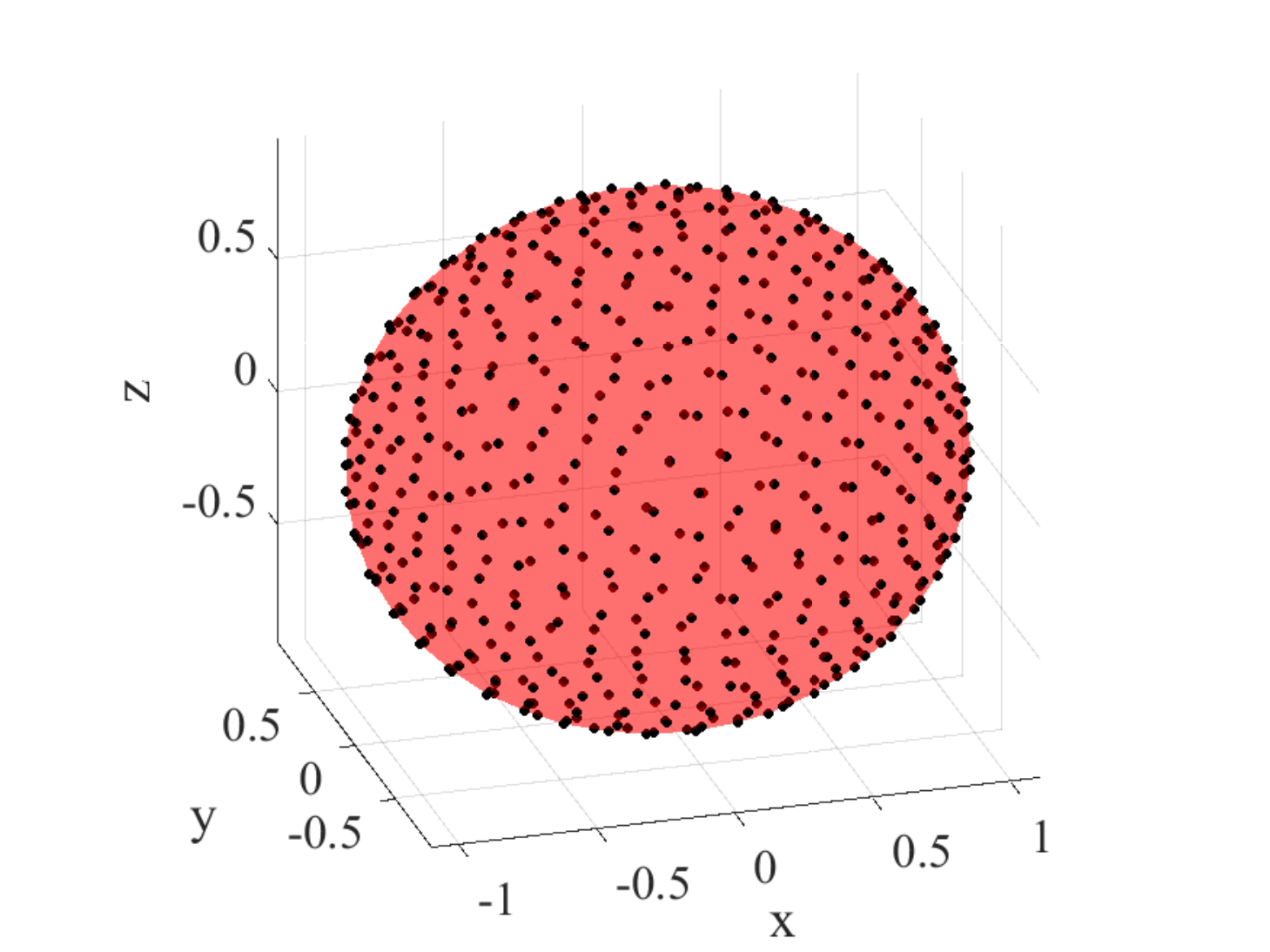}}
\subfigure[Perturbed ellipsoid]{\label{fig:shapeb}\includegraphics[width=75mm]{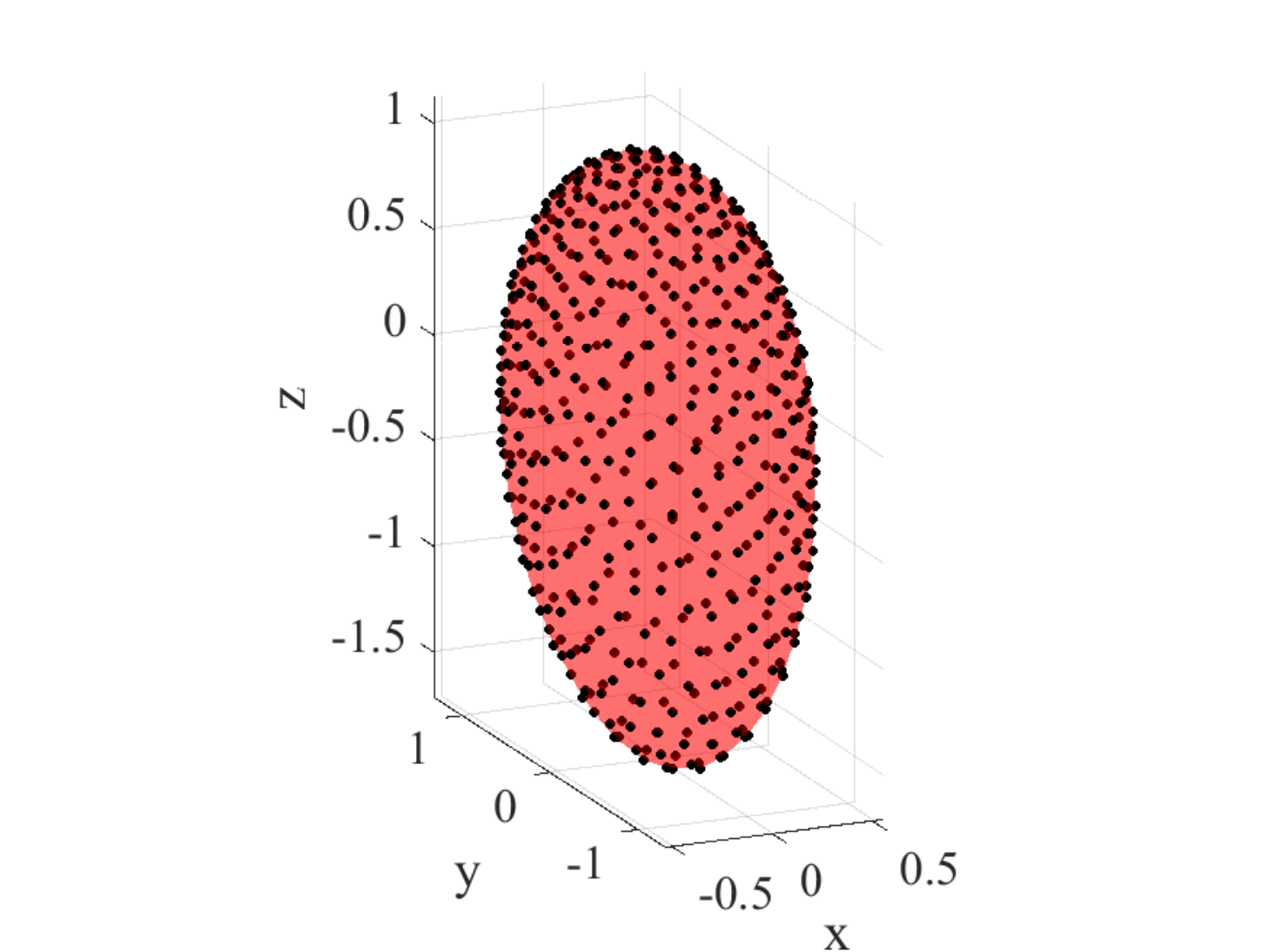}}
\caption{The shapes for study. 529 evaluation points are shown. The distribution of other sets of evaluation points is similar.}
\label{fig:shapes}
\end{figure}
We use parameters $a=1.1$, and  $b=c=\fraco{1}{\sqrt{1.1}}$ in Eq.\ \eqref{eq:obj1} and $a=0.1,\, b=c=0.2$, and $B=0.25$ in Eq.\ \eqref{eq:obj2}. In both cases, the reference configuration is the unit sphere, $\vec{Z}=(\cos{\lambda}\cos{\theta}, \sin{\lambda}\cos{\theta}, \sin{\theta})$, and the parameter $V_f$ in Eq.\ \eqref{eq:obj2} is set so that the volume of the perturbed ellipsoid is the same as that of the unit sphere (in our case, $V_f \approx 5.14$). For simplicity, we also let $\sigma=K=A=G_s=1$ in the force and energy functionals in Eqs.\ \eqref{eq:nhforce} and \eqref{eq:stforce}. The shapes are shown in Fig.\ \ref{fig:shapes} along with the positions of the set of 529 evaluation points (sets with a higher number of evaluation points are similarly distributed). 

\begin{figure}[!tbp]
\centering  
\includegraphics[width=90mm]{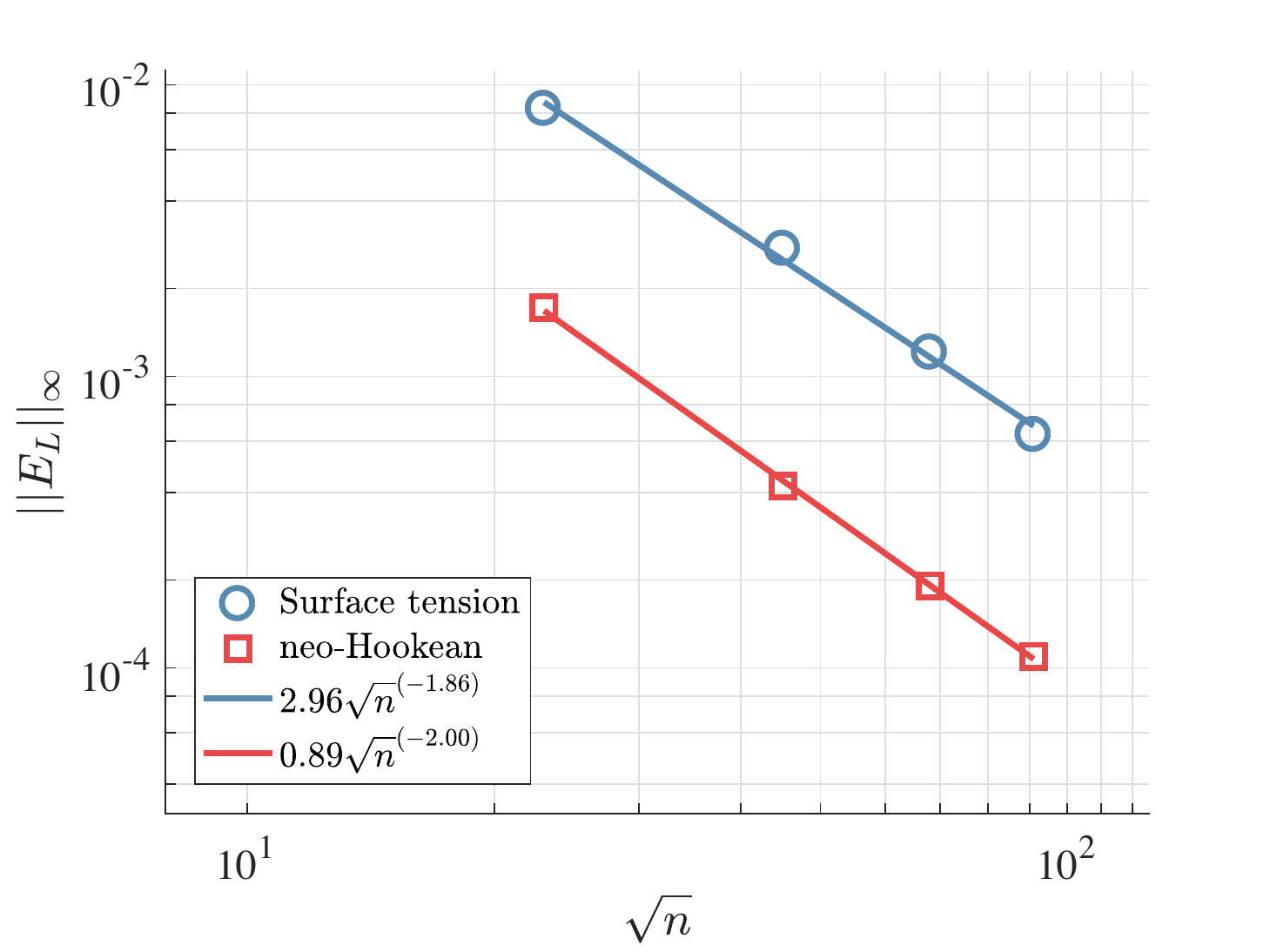}
\caption{Errors in surface tension and neo-Hookean forces on an ellipsoid. $\norm{\vec{E_L}}$, the $\ell_\infty$ error in LDSM forces, is shown for the cases of 529, 2025, 4624, and 8281 points along with a power-law fit of the data. Because the spherical harmonic interpolant from Eq.\ \eqref{eq:shintetp} is exact for an ellipsoid, forces computed from SHVD are exact and $\norm{\vec{E_S}}=0$.}
\label{fig:ellipseerror}
\end{figure} 

We fix a number of evaluation points so that the quadrature error given by the weights in Eq.\ \eqref{eq:wtseq} is on the order $10^{-5}$ or less. The quadrature error is measured by the value of the discrete force integral since the exact integral of the translation invariant continuous force function evaluates to zero. We then compute the exact solution for the neo-Hookean or surface tension force densities at the $n$ evaluation points using Eq.\ \eqref{eq:vdforce} to take the variational derivative of \textit{the exact} mappings $\vec{X}$ from Eqs.\ \eqref{eq:obj1} or \eqref{eq:obj2}.

Next, we construct a spherical harmonic interpolant for the surface of the form in Eq.\ \eqref{eq:shintetp} using some number $m$ of interpolation points. We then use this interpolant to compute the force densities in Eq.\ \eqref{eq:vdforce}. We multiply both the exact force density and the force density from SHVD by the weights in Eq.\ \eqref{eq:wtseq} to get a force. In this way, we isolate the error from representing the surface with a spherical harmonic interpolant. Finally, we compute the LDSM forces by evaluating Eq.\ \eqref{eq:faiforce} with the triangle vertex locations given by Eq.\ \eqref{eq:obj1} or \eqref{eq:obj2}. 

We denote the exact force (not force density) as $\vec{\widehat{F}}_E$, the force from SHVD as $\vec{\widehat{F}}_S$ and the force from LDSM as $\vec{\widehat{F}}_L$. We quantify the error between the models via the discrete $\ell_\infty$ norm. The norm of the errors from SHVD and LDSM, denoted by $\vec{E_S}$ and $\vec{E_L}$, respectively, are given by
\begin{gather}
\label{eq:errorshvd}
\norm{\vec{E_S}}_\infty = \max_i{\norm{\vec{\widehat{F}}_S - \vec{\widehat{F}}_E}_{\infty}}
\end{gather}
and
\begin{gather}
\label{eq:errorldsm}
\norm{\vec{E_L}}_{\infty} = \max_i{\norm{\vec{\widehat{F}}_L - \vec{\widehat{F}}_E}_{\infty}}.
\end{gather}
The error in force is computed as the maximum over the $n$ evaluation points, where $i$ is the index that runs from 1 to $n$. 
In all convergence plots, we plot the error vs. $\sqrt{m}$ or $\sqrt{n}$. Since $(dA_0)_i \approx 4\pi/n$, the number of points is proportional to the square of the length scale. We plot convergence in the length scale by taking the square root of the number of points. 

\begin{figure}[!t]
\centering     
\subfigure[LDSM]{\label{fig:LDSMpe}\includegraphics[width=75mm]{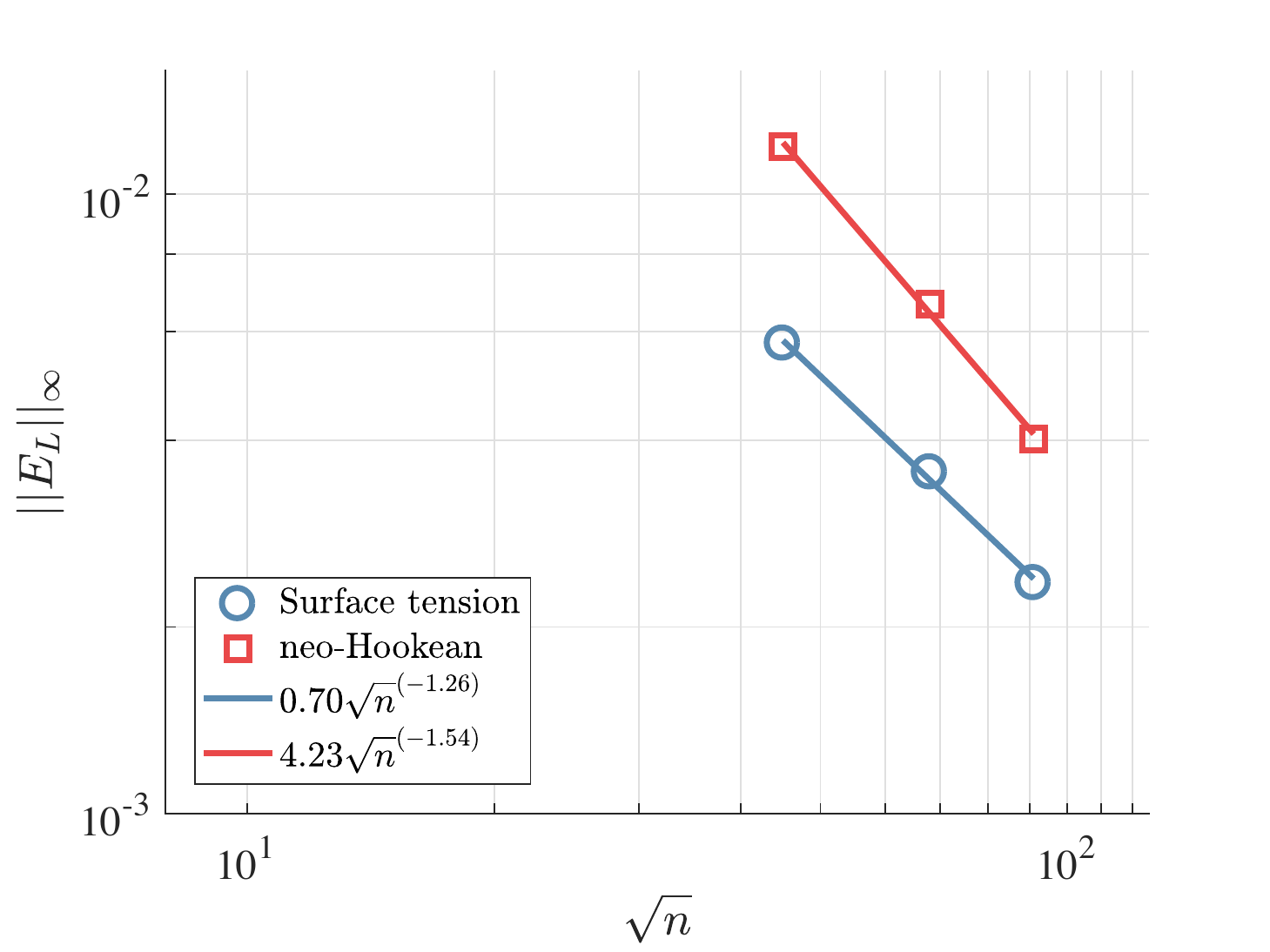}}
\subfigure[SHVD]{\label{fig:SHVDpe}\includegraphics[width=75mm]{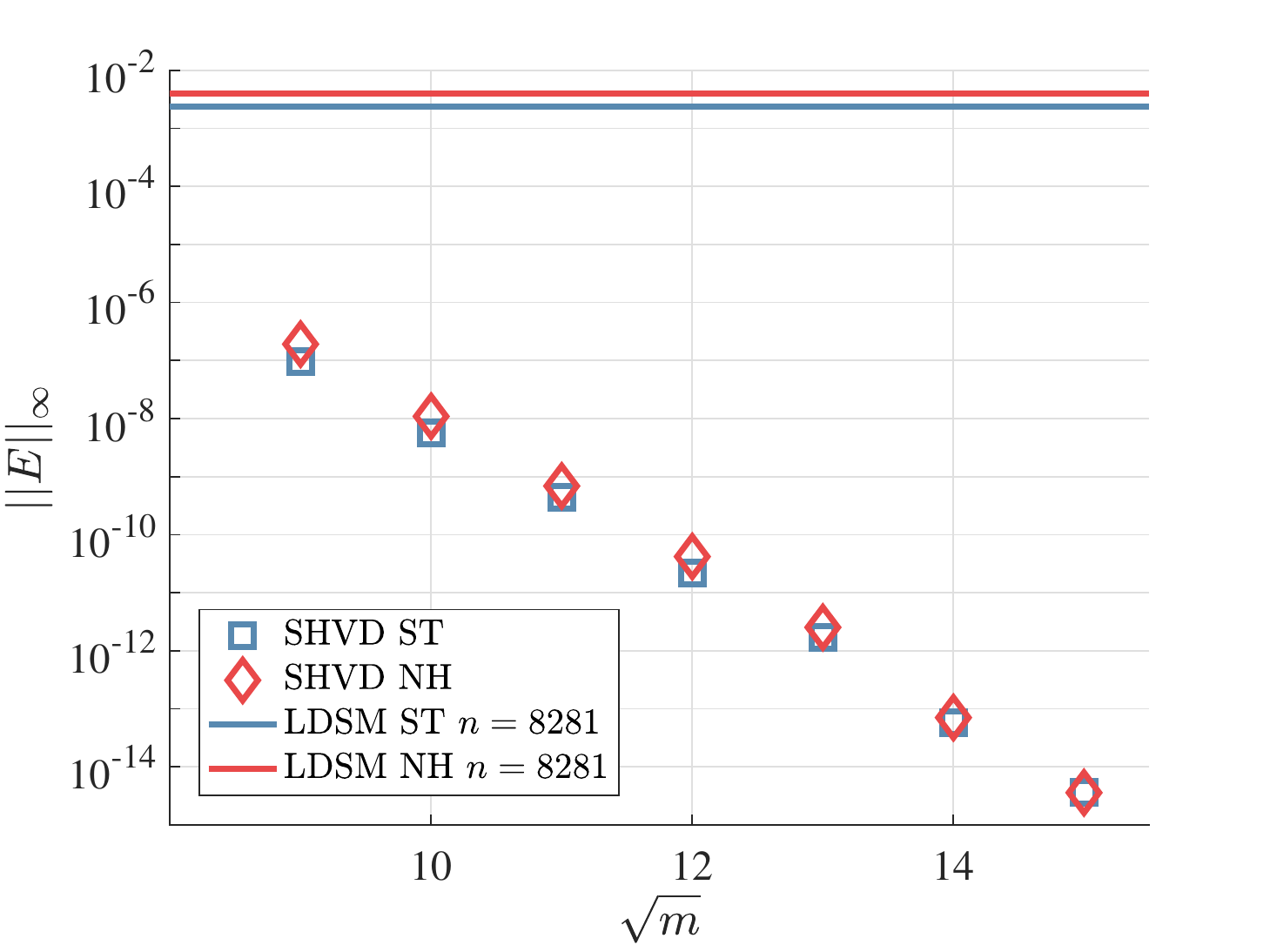}}
\caption{Errors in forces for the perturbed ellipsoid (Eq.\ \eqref{eq:obj2}) for the two models. The discrete $\ell_\infty$ norms of $\vec{E_L}$ and $\vec{E_s}$ given by Eqs. \eqref{eq:errorldsm} and \eqref{eq:errorshvd} are graphed. In (a), the error $\norm{\vec{E_L}}$ is plotted against the number of evaluation points $n$ for both surface tension (blue circles) and neo-Hookean (red squares) forces. In (b), the error in SHVD $\norm{\vec{E_S}}$ is plotted against the number of interpolation points $m$ for both surface tension (blue squares) and neo-Hookean (red diamonds) forces. We observe spectral convergence in the number of interpolation points $m$. Forces computed with SHVD with as few as $m=81$ interpolation points are
several orders of magnitude more accurate than those computed from LDSM with $n=8281$ evaluation points.}
\label{fig:forceerrorsPE}
\end{figure}

Because the same area weights are used for the calculation of ${\vec{\widehat{F}}_E}$ and ${\vec{\widehat{F}}_S}$, 
the measured error in Eq.\ \eqref{eq:errorshvd} is independent of the quadrature weights. Since LDSM uses different area weights, the measured error in Eq.\ \eqref{eq:errorldsm} is affected by both the quadrature error in $\vec{\widehat{F}}_E$ and the discretization error in the LDSM force calculation. We define $\vec{\widehat{F}}_E$ to be the exact solution if the quadrature error (i.e. the distance of the force integral from zero) is negligible relative to the LDSM surface discretization error. We only consider surfaces where the quadrature error is $\mathcal{O}(10^{-5})$ or less. In our numerical tests, the LDSM  error graphed in Figs.\ \ref{fig:ellipseerror} and \ref{fig:LDSMpe} is larger than this value so that quadrature error is negligible. 


Fig.\ \ref{fig:ellipseerror} shows the errors in Eqs.\ \eqref{eq:errorshvd} and \eqref{eq:errorldsm} for the forces on the ellipsoid as a function of the number of evaluation points, for $n=$ 529, 2025, 4624, and 8281. We note that for an ellipsoid in Eq.\ \eqref{eq:obj1}, the spherical harmonic interpolant with as little as 4 interpolation points gives the same mapping as Eq.\ \eqref{eq:obj1}. Therefore, the forces from SHVD are exact for an ellipsoid and are not included in Fig.\ \ref{fig:ellipseerror}. 
Consequently, we show in Fig.\ \ref{fig:ellipseerror} only the errors in LDSM as a function of the number of evaluation points for both surface tension and neo-Hookean forces. A piecewise planar representation of a surface converges to a smooth surface with second order accuracy, so we expect at most second order accurate LDSM forces. We observe a second order convergence rate in the forces with respect to the length dimension (see power function fit of data in Fig.\ \ref{fig:ellipseerror}). We observe a similar convergence rate for computing surface tension on a sphere with LDSM (data not shown).

Our study of the perturbed ellipsoid shows that forces from SHVD converge to the exact solution with spectral accuracy in contrast to LDSM. Fig.\ \ref{fig:forceerrorsPE} shows the errors in surface tension and neo-Hookean forces for LDSM with varying numbers of evaluation points $n$ (Fig.\ \ref{fig:LDSMpe}) and SHVD with varying numbers of interpolation points $m$ (Fig.\ \ref{fig:SHVDpe}). A convergence rate between first and second order is observed for LDSM. In contrast, the force from SHVD converges spectrally as the number of interpolation points $m$  increases regardless of the number of evaluation points $n$ \cite{shankar}.
The force estimate from SHVD is orders of magnitude smaller than LDSM when $m \geq 81$. This stems from the fact that SHVD relies on a \textit{continuous} surface representation, and is therefore able to effectively capture the geometrical information of the surface necessary to accurately compute forces. 

We comment that the error in forces from SHVD is from both interpolation error and aliasing, a well known phenomenon in which coefficients of higher order spherical harmonics are folded into lower order coefficients computed in Eq.\ \eqref{eq:systeminterp} \cite{alias}. Here the analytical expressions for force are often not in the space of spherical harmonics of degree at most $N$, which makes them subject to aliasing error. Previous studies, among them \cite{veerapaneni2011fast,zhao2010spectral}, have accounted for this via upsampling of points for force calculations. We take a different approach by increasing the number of interpolation points until the error becomes close to machine epsilon. Fig.\ \ref{fig:SHVDpe} shows the error in the force calculation for the perturbed ellipse to be on the order $10^{-15}$ with 225 interpolation points used.

\subsection{Dynamic test}
\label{sec:dynamics}
%
%
\begin{figure}
\centering  
\includegraphics[width=90mm]{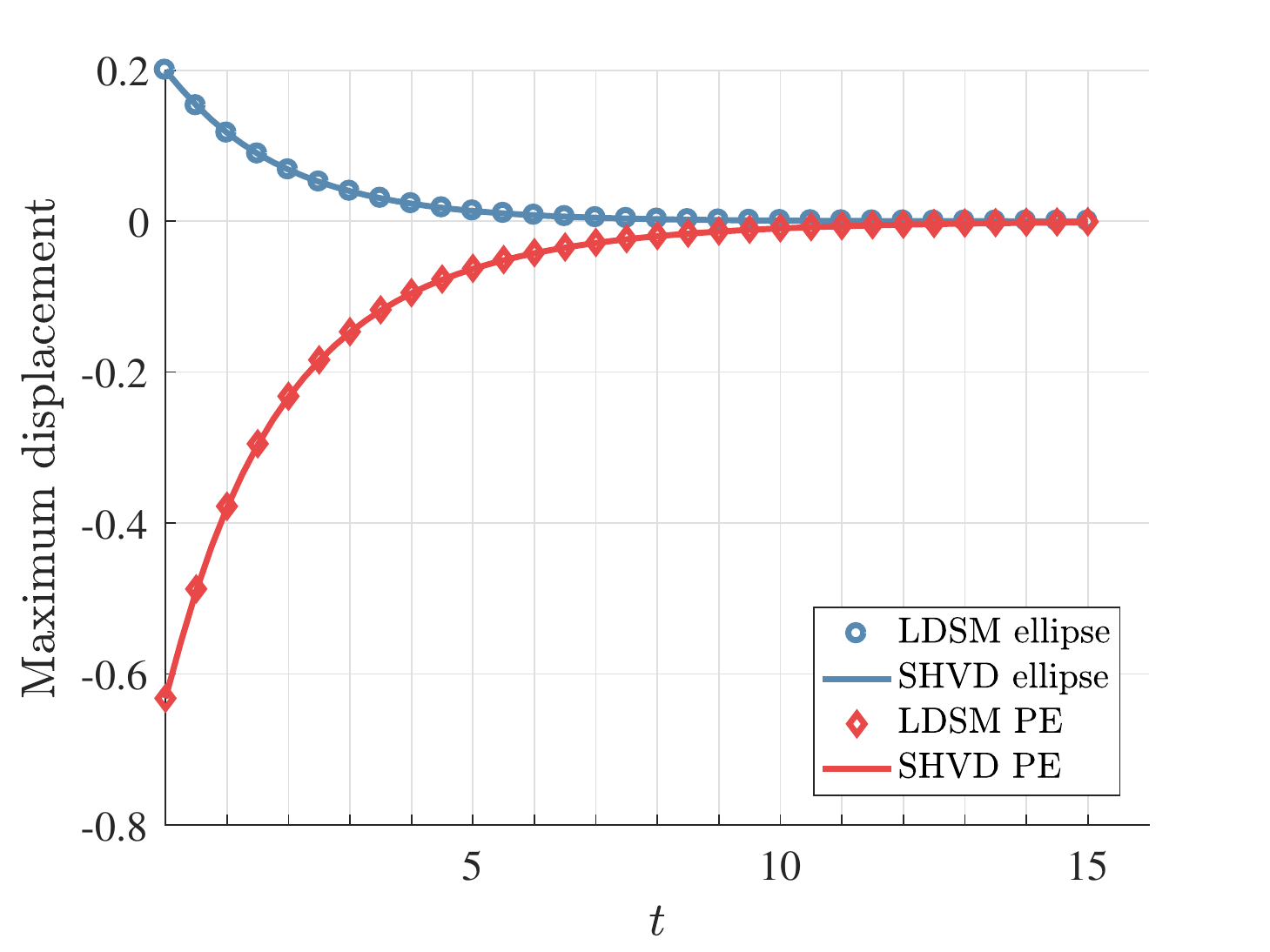}
\caption{Maximum displacement over time for dynamic simulations of the test objects. The displacement is the maximum displacement from the center of mass minus the radius in the reference configuration. Despite the inaccuracy of LDSM in computing the forces, the full IB simulation gives similar results to one that uses SHVD for force computations.}
\label{fig:Balloons}
\end{figure}

We next perform full immersed boundary simulations on the test objects in Eqs.\ \eqref{eq:obj1} and \eqref{eq:obj2} with parameters $a=1.2, \, b=c=1/\sqrt{1.2}$, $B=0.25$, $V_f \approx 5.14$, and $\sigma=K=A=G_s=1$. 
Each object relaxes to its equilibrium configuration  over time (in this case, the unit sphere). The relaxation timescale is dictated by the values of elastic moduli and the viscosity of the fluid, which we set to $\mu=1$ for simplicity.
We use 8281 evaluation points ($n = 8281$) for both objects and 225 interpolation points in the SHVD model ($m=225$). 
From Fig.\ \ref{fig:ellipseerror}, this corresponds to an initial error in forces on the order of $10^{-3}$ for LDSM and $10^{-15}$ (machine epsilon) for SHVD. We use a total force that is the sum of surface tension and neo-Hookean forces (i.e. $\vec{\widehat{F}}_{TOT}$ = $\vec{\widehat{F}}_{NH} + \vec{\widehat{F}}_{ST}$). The Eulerian grid runs from $a=-2$ to $b=2$ in the $x,y$, and $z$ directions with a grid size of $\eta=32$ (32 Eulerian points in each direction). For each object, we track the maximum or minimum displacement from the center of mass as a function of time. Fig.\ \ref{fig:Balloons} shows that, despite the large difference in accuracy between LDSM and SHVD force computations, both models give similar results for displacement over time due to the first order time update and the discrete delta function in the spreading and interpolation operators.
\begin{figure}
\centering  
\includegraphics[width=90mm]{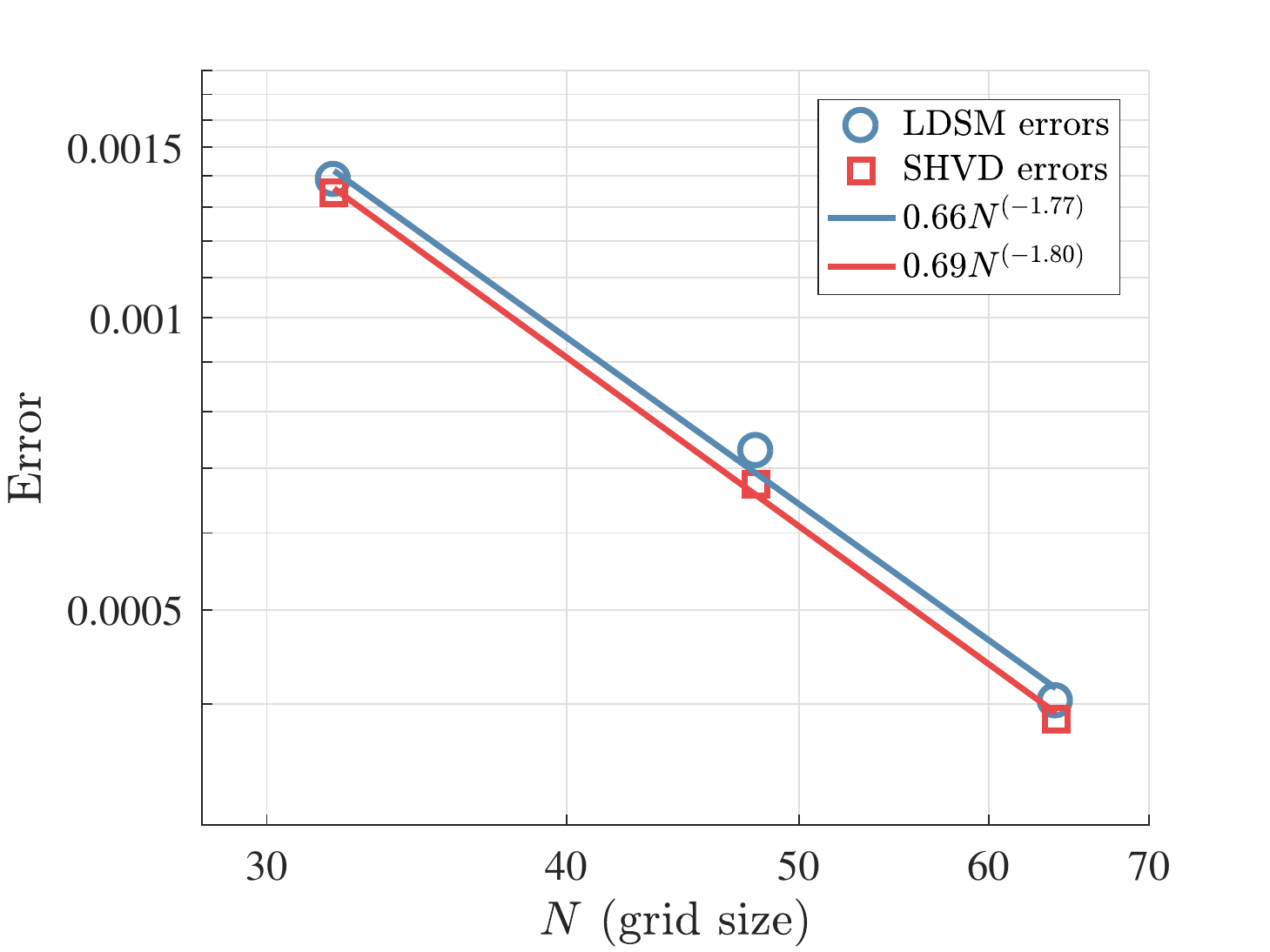}
\caption{Convergence tests on a dynamic simulation of an ellipsoid for both LDSM and SHVD. The error is measured as the change in displacement of the point with largest $x$ coordinate relative to the center of mass, $d_{max}$. The error is the difference in $d_{max}$ with respect to a grid 1.5 times as refined. Lines show power law function fits to the data which demonstrate the approximately first order accuracy of both methods.}
\label{fig:balloonconv}
\end{figure}

We study the convergence of each method by examining the state of the ellipsoid in Eq.\ \eqref{eq:obj1} and Fig.\ \ref{fig:shapea} at $t=3$ seconds. The initial surface configuration is the ellipsoid with $a=1.2$, and $b=c=1/\sqrt{1.2}$. We consider a fixed Eulerian grid of size $\eta^3$ Eulerian points. We then choose a Lagrangian grid with $n \approx 2\eta^2$ points. We let $\Delta t = 1/(2\eta)$ and measure the change in distance, denoted $d_{max}$, of the rightmost point (the point with largest $x$ coordinate) minus the center of mass (i.e. it is the amount the rightmost point has moved in the $x$ direction after 3 seconds after subtracting the center of mass). The error is the difference in $d_{max}$ with respect to a grid for the same method that is 1.5 times as refined. For example, the error value for a 32 point LDSM grid is the change in $x$-displacement of the rightmost point from the center of mass on the 32 point grid minus the corresponding value on a 48 point grid, both using LDSM. 

We consider grid sizes $\eta$, of $32$, $48$, and $64$ corresponding to $n=2025,\, 4624$, and $8281$ evaluation points. For SHVD, we hold the number of interpolation points $m$ constant at $64$.
Fig.\ \ref{fig:balloonconv} shows the results of the convergence test for the different grid sizes. Both methods display at least first order convergence in space and time. This is expected for the IB method, which as discussed in \cite{lai2000immersed,griffith2005order},  is limited to first order accuracy in the case of an infinitely thin elastic membrane. We comment on this further in Section \ref{sec:discuss}.

In addition to studying the convergence of each method under these circumstances, we also compare the computational costs. In a similar procedure to the convergence tests, we track the ratio of SHVD to LDSM execution times (i.e. time for SHVD / time for LDSM) for a 15 second ellipsoidal contraction with the same parameters as the convergence study. Table \ref{tab:times} summarizes the results. For $m$ two or more orders of magnitude less than $n$, the wall clock time of our method is within 20\% of LDSM. Observe that for larger $m$ values, the $\mathcal{O}(mn)$ operations in our method's force computation increase the computational time farther beyond that of LDSM. We will see in Sections \ref{sec:RBCs} and \ref{sec:bleb} how this behavior impacts the use of the model in our applications. For $m \ll n$ (the scenario under which our method is intended for use), the algorithms have comparable execution times. 
\begin{table}
\centering
\begin{tabular}{|cccc|} 
\hline
$\eta$ & $n$ & $t$-ratio for $m=64$ & $t$-ratio for $m=225$ \\
\hline
32 & 2025 & 1.15 & 1.73 \\
48 & 4624 & 1.08 & 1.46 \\
64 & 8281 & 0.99 & 1.30 \\
72 & 10404 & 1.00 & 1.32 \\
\hline
\end{tabular}
\caption{Ratio of computational times for 15 second ellipsoidal contraction.}
\label{tab:times}
\end{table}

%
\section {Applications}
\label{sec:apps}
In this section, we consider two different applications involving cellular dynamics: red blood cells in capillary flow and a model of bleb expansion. The former application involves large and smooth deformations of the entire cell surface due to the background flow in a narrow tube. The membrane deformation in the model of cellular blebbing occurs because of a discontinuous applied force on a small region of the membrane surface. For both models, we compare simulation results from our method to LDSM. The results from the model simulations highlight the versatility of our approach.

\subsection{Red blood cells in capillary flow}
\label{sec:RBCs}
%
%
Red blood cells (RBCs) passing through microcapillaries is a well-studied biophysical problem that is important for understanding the rheological properties of both large and micro scale blood flow as well as transport via blood flow \cite{pozrikidis2005axisymmetric}. It has been observed experimentally that RBCs, which are nearly incompressible, deform into parachute-like shapes when flowing through capillaries, with the leading, or front, edge protruding under the influence of Poiseuille flow within the capillary \cite{tsukada2001direct}.
For smaller tube diameters, it has been shown that RBCs align themselves in a single file configuration \cite{tsukada2001direct}. Our approach is therefore to model a single RBC flowing through a long capillary with SHVD. Because the velocity and lengthscales are small, the Reynolds number for flow in the microcirculation is on the order $< 10^{-2}$, and the Stokes equations, Eqs.\ \eqref{eq:Stokeseqnsone}-\eqref{eq:Stokeseqnstwo}, govern the flow profile. We also simplify the problem by assuming the viscosity of the internal hemoglobin solution to be equal to that of the surrounding plasma \cite{pozrikidis2005axisymmetric}. 
\subsubsection{Model of the RBC and capillary}
The reference configuration of a red blood cell that we use is taken from \cite{omori2012tension}. Let $(x_s, y_s, z_s) = (\cos{\lambda} \cos{\theta}, \sin{\lambda} \cos{\theta}, \sin{\theta})$ denote the mapping from $(\lambda, \theta) \rightarrow \mathbb{R}^3$ for the unit sphere. The RBC reference configuration is given by
\begin{equation}
\label{eq:zrbc}
\vec{Z_{RBC}} = R_{RBC}\begin{pmatrix}
\fraco{1}{2}x_s\left(0.21 + 2\left(y_s^2 + z_s^2\right) - 1.12\left(y_s^2+z_s^2\right)^2\right)\\[4 pt]
y_s\\[4 pt]
z_s
\end{pmatrix}.
\end{equation}
The mean radius of the RBC is $R_{RBC}=3.91$ $\mu$m \cite{pozrikidis2005axisymmetric}. Elastic forces on the RBC come from resistance to stretching, shear, and bending. We use the neo-Hookean energy functional in Eq.\ \eqref{eq:nhforce} with $G_s = 2.5 \times 10^{-3}$ dyn/cm $=2.5$ pN/$\mu$m and $A=50$ to model resistance to shear and stretching. As discussed in \cite{fai, pozrikidis2005axisymmetric}, $A$ is a computational parameter that is used to enforce the incompressibility of the RBC membrane while ensuring stability at large enough timesteps. \cite{fai} reported using $A$ in the range from $8$ to $400$, and we found our value of $A=50$ to give area deformations less than 1\%. 
The derivatives of the reference configuration required for the neo-Hookean force calculation in Eq.\ \eqref{eq:vdforce} are computed analytically from the $\vec{Z_{RBC}}$ mapping. In addition, the area weights required to integrate the force over the RBC reference configuration are given by Eq.\ \eqref{eq:arbitwts}, where $\vec{Z}_\mathbb{L}=\vec{Z_{RBC}}$. 
For bending, we use the energy functional in Eq.\ \eqref{eq:bendenergy} with a typical value $k_{bend} = 1 \times 10^{-12}$ erg $=0.1$ pN $\cdot$ $\mu$m. 

The capillary is modeled as a thin elastic cylindrical membrane of radius 5 $\mu$m that is discretized and tethered in place. Letting $\vec{Z}_{cap}$ and $\vec{X}_{cap}$ be the reference and deformed configurations of the capillary, the Lagrangian \textit{force} on a capillary node is given by
\begin{equation}
\label{eq:tetherF}
\vec{\widehat{F}}_{teth}(\vec{q}) = -k_{teth}\left(\vec{X}_{cap}(\vec{q}) - \vec{Z}_{cap}(\vec{q})\right), 
\end{equation}
where $\vec{q}$ refers to the Lagrangian coordinate of the node and $k_{teth} = 2.45 \times 10^{-3}$ dyn/cm $= 2.45$ pN/$\mu$m. 
As discussed in \cite{teran2009tether}, tether force functions such as Eq.\ \eqref{eq:tetherF} are not translation invariant, and are thus not guaranteed to integrate to zero in the discrete sense. Because the solution of the Stokes equations, Eqs.\ \eqref{eq:Stokeseqnsone}-\eqref{eq:Stokeseqnstwo}, on a periodic domain is non-unique, a constant velocity $\vec{u_c}(t)$ can be added to the structure velocities obtained from interpolation at each time step. Referring to Section \ref{sec:timeupdate}, the time-stepping procedure is modified so that in the sixth step, the cylinder boundary is updated to an intermediate position
\begin{equation}
\hat{\vec{X}}_{cap}^{k+1} = \vec{X}^{k}_{cap} + \Delta t \, \vec{U}^{k+1}.
\end{equation}
Then $\hat{\vec{X}}_{cap}^{k+1}$ is used to calculate a constant velocity that is added throughout the entire domain (also to the interpolated velocity of the RBC) to ensure that the tether forces sum to zero at the next time step. The constant velocity is given by
\begin{equation}
\label{eq:ucorrect}
\vec{u_c}(t) = \fraco{1}{A_{cap} \Delta t} \int_{cap} \left(\vec{Z}_{cap}(\vec{q}) - \hat{\vec{X}}_{cap}^{k+1}(\vec{q})\right) \, d\vec{q}, 
\end{equation}
where $A_{cap}$ is the (lateral) surface area of the capillary. A derivation of Eq.\ \eqref{eq:ucorrect} can be found in \cite{teran2009tether}. The positions of the Lagrangian structures (capillary and RBC) are then updated by $\vec{X}^{k+1} = \hat{\vec{X}}^{k+1} + \vec{u_c}(t) \Delta t$ to ensure the tether forces in Eq.\ \eqref{eq:tetherF} integrate to zero at the next time step. 

A background force of $0.08$ $\textrm{pN}/\mu \textrm{m}^3$ $= 8.0 \times 10^3$ $\textrm{dyn/cm}^3$ is specified inside the capillary to establish a parabolic flow profile as in \cite{fai}. Given the viscosity of blood plasma as 1.2 cP = $1.2 \times 10^{-3} \textrm{ pN} \cdot \textrm{s} / \mu \textrm{m}^2$, the resulting maximum flow velocity is 750 $\mu$m/s inside the capillary as in \cite{fai,tsukada2001direct}. We impose a uniform force in the opposite direction outside of the capillary to ensure the integral of the background force over the Eulerian grid is zero. 
We use a $24$ $\mu$m $\times$ $24$ $\mu$m $\times$ $24$ $\mu$m fluid grid. The long axis of the capillary is the $x$-axis, and the capillary runs from $x=-8$ to $x=8$ $\mu$m. We center the cell at $(-8, 0, 0)$ initially and simulate its deformation until it reaches a steady state shape (around $t=0.05$ s $=500 \Delta t$). 

\subsubsection{Results}
As discussed in Section \ref{sec:brforce}, there is no standard way to model bending energy in LDSM because the piecewise linear triangles have zero curvature. In addition, the bending force model from \cite{fai} is not guaranteed to integrate to zero in the discrete sense, meaning it cannot be used in a zero Reynolds number regime with periodic boundary conditions. Thus in order to fairly compare our model to LDSM, we begin by examining the case of $k_{bend}=0$. For all simulations, we use $n=8281$ evaluation points and vary the number of interpolation points $m$ in our method.  
\begin{figure}[!htbp]
\begin{center}
\includegraphics[width = 0.9\textwidth]{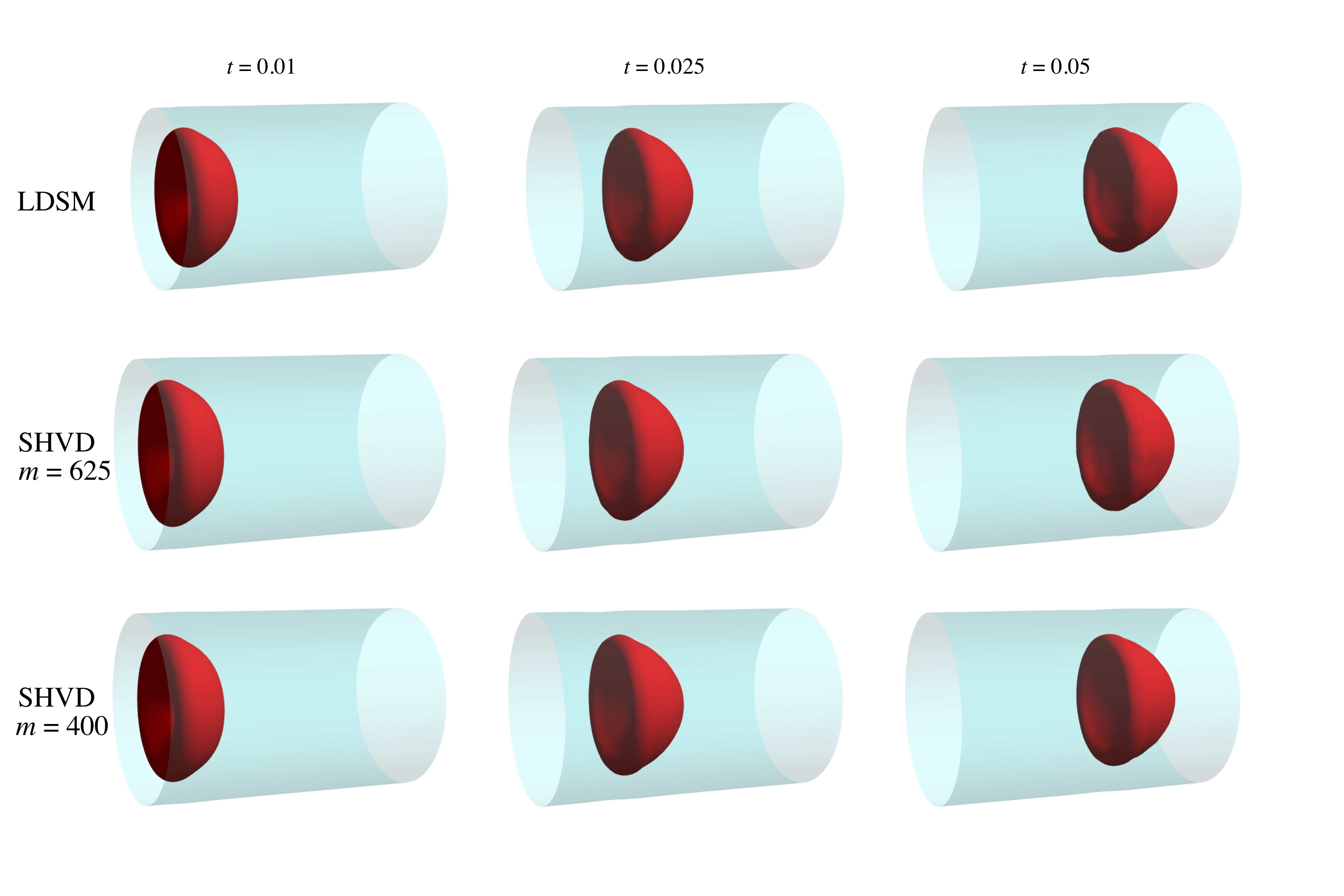}
\caption{The cell membrane and capillary at several time values using LDSM (top) and SHVD with $m=625$ (middle) and $m=400$ (bottom) interpolation points. All methods used $n=8281$ evaluation points.}
\label{fig:rbcincylinder}
\end{center}
\end{figure}

Fig.\ \ref{fig:rbcincylinder} shows the deformation of the RBC from its biconcave initial shape into the parachute-like shapes observed in experimental capillary flow \cite{tsukada2001direct} for LDSM (top), and SHVD with $m=625$ (middle) and $m=400$ (bottom). 
We refer to the shape at $t=0.05$ as the steady state shape of the RBC because it does not change significantly after this time value.
While the steady state shapes are similar, a closer look reveals bumps in the RBCs modeled with LDSM and SHVD with $m=625$. We take these asymmetric bumps to be non-physical, numerical error since the imposed flow and cylinder are symmetric. These bumps are distinct from the larger wavelength, axisymmetric, wrinkles observed in all of our simulations that exclude bending rigidity. We note that wrinkles have been observed in previous studies of shear flow at low bending rigidities \cite{fai,yazdani2011phase}. 
When the number of interpolation points is greater than 625, we observe even more bumps on the cell (data not shown), indicating instability due to higher frequency modes in the interpolant. For the remainder of this section, we take $m=400$ to obtain smooth surfaces similar to previous computational studies of RBCs \cite{fai,veerapaneni2011fast}.

\begin{figure}[!htbp]
\begin{center}
\includegraphics[width = 0.75\textwidth]{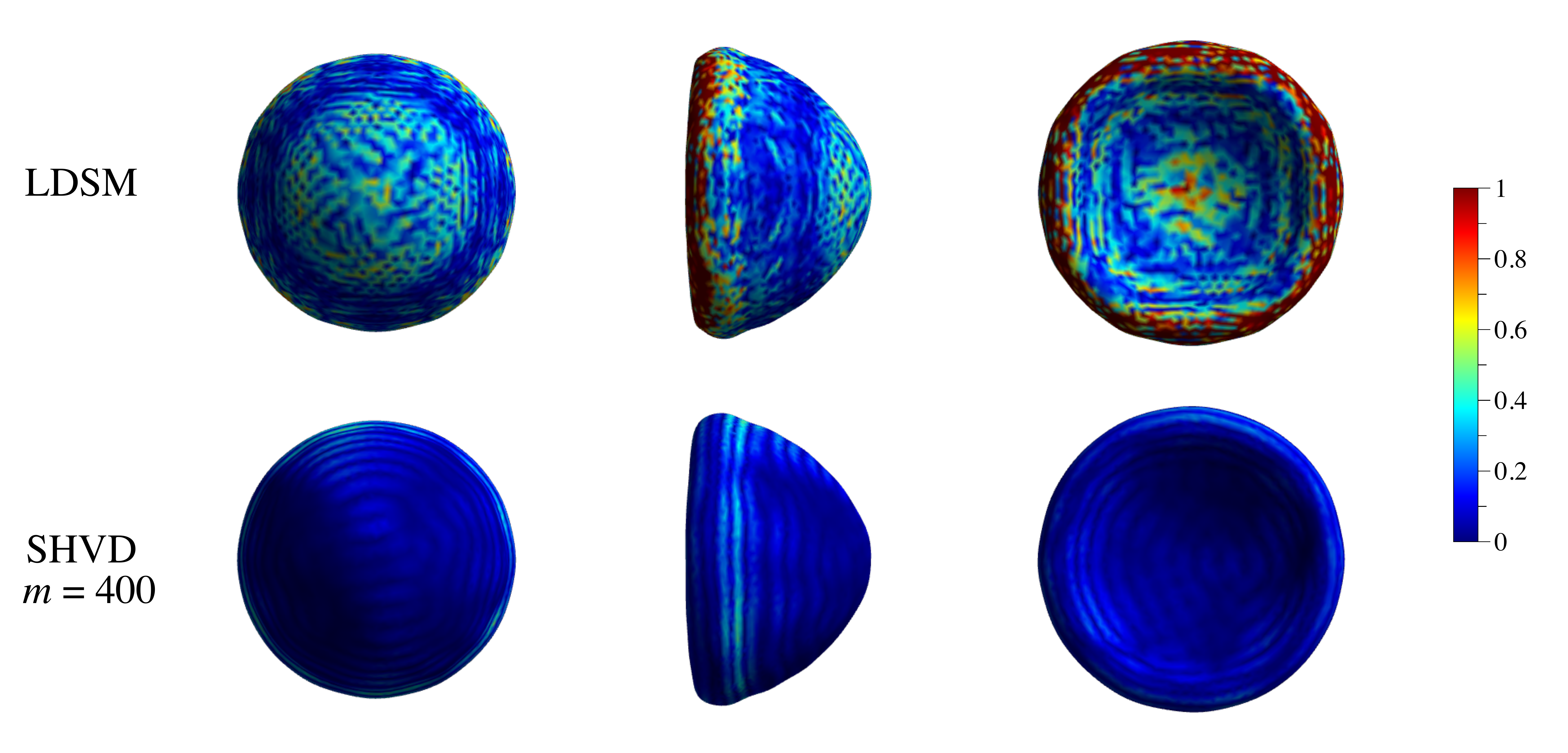}
\caption{Magnitude of the elastic forces on the cell at the front, side, and back views of the cell at $t=0.05$ from simulations with LDSM (top) and SHVD with $m=400$ interpolation points (bottom). $n=8281$ evaluations points were used for both simulations.}
\label{fig:rbcforcecompare}
\end{center}
\end{figure}

Fig. \ref{fig:rbcforcecompare} shows a comparison of the magnitude of the RBC elastic forces from simulations using LDSM and SHVD at $t=0.05$. Although the overall accuracy in space and time of both methods is first order (see Section \ref{sec:dynamics}), the forces are much smoother with SHVD. Additionally, the magnitude of the forces is much larger in computations with LDSM than SHVD at the back end of the cell when its shape changes from concave to convex (and the mean curvature changes sign). The side views of the steady state shape (middle, Fig.\ \ref{fig:rbcforcecompare}) show wrinkles that occur in shapes from both simulations when the bending forces are excluded. 

\begin{figure}[!htbp]
\begin{center}
\includegraphics[width = 0.9\textwidth]{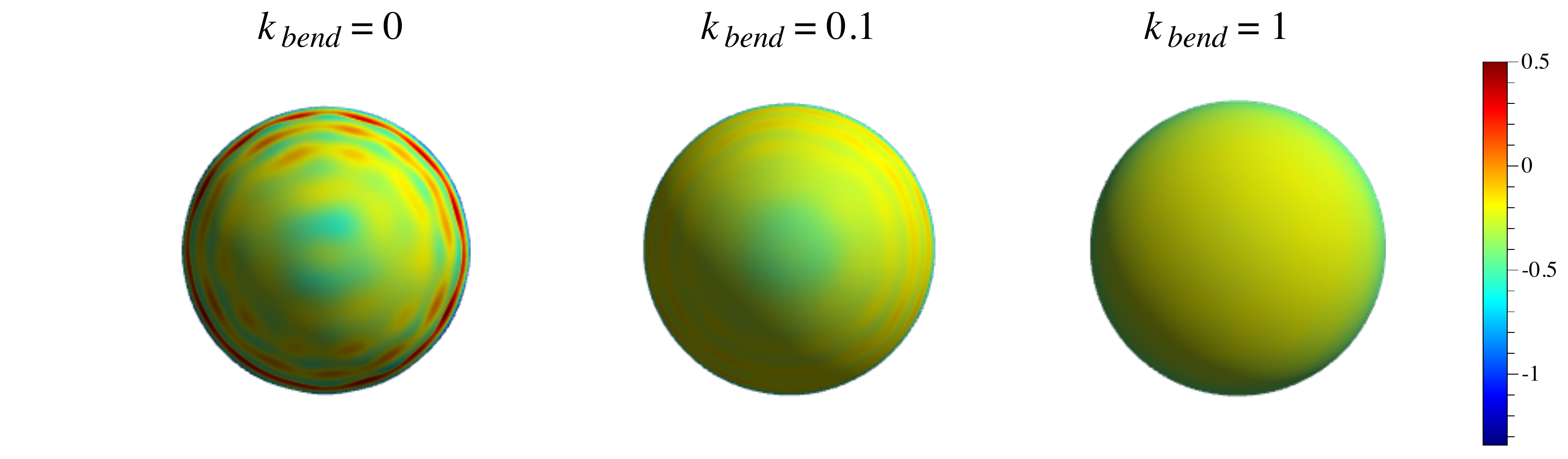}
\caption{Mean curvature of the cell viewed from the front for increasing values of $k_{bend}$. The cell was simulated with SHVD ($m=400$ interpolation and $n=8281$ evaluation points).}
\label{fig:curvcompare}
\end{center}
\end{figure}

Surface bumps and wrinkles of the RBC are visible in Fig.\ \ref{fig:rbcforcecompare} when $k_{bend}=0$. Fig.\ \ref{fig:curvcompare} shows a head-on (facing the cell) view of the steady state cell shapes (at $t=0.05$) for different values of bending rigidity. 
When $k_{bend}=0$, the mean curvature $H$ (defined in Eq.\ \eqref{eq:meancurve}) oscillates between positive and negative because of the wrinkling 
(Fig.\ \ref{fig:curvcompare}, $k_{bend}=0$).
When $k_{bend}=0.1$ pN $\cdot$ $\mu$m, surface wrinkles are significantly reduced, the curvature remains negative on the top face, and a smooth parachute shape is obtained. We note this is the same value of $k_{bend}$ used in \cite{fai}. The curvature remains negative on the top face when $k_{bend}=1$.

At steady state, the front of cells with higher bending rigidities more closely resembles a sphere. We quantify this by examining the variability in mean curvature on the front edge of the cell as a function of the bending rigidity in Fig. \ref{fig:curvcompare}. It can be seen that as $k_{bend}$ increases from 0 to 1 pN $\cdot$ $\mu$m, the front of the cell more closely resembles a sphere with uniform curvature, as is expected for large bending rigidities. We note the smoothness of the curvature values in our method, which clearly shows the trend towards a more spherical shape for larger bending rigidities.

\begin{figure}[!htbp]
\begin{center}
\includegraphics[width = 0.9\textwidth]{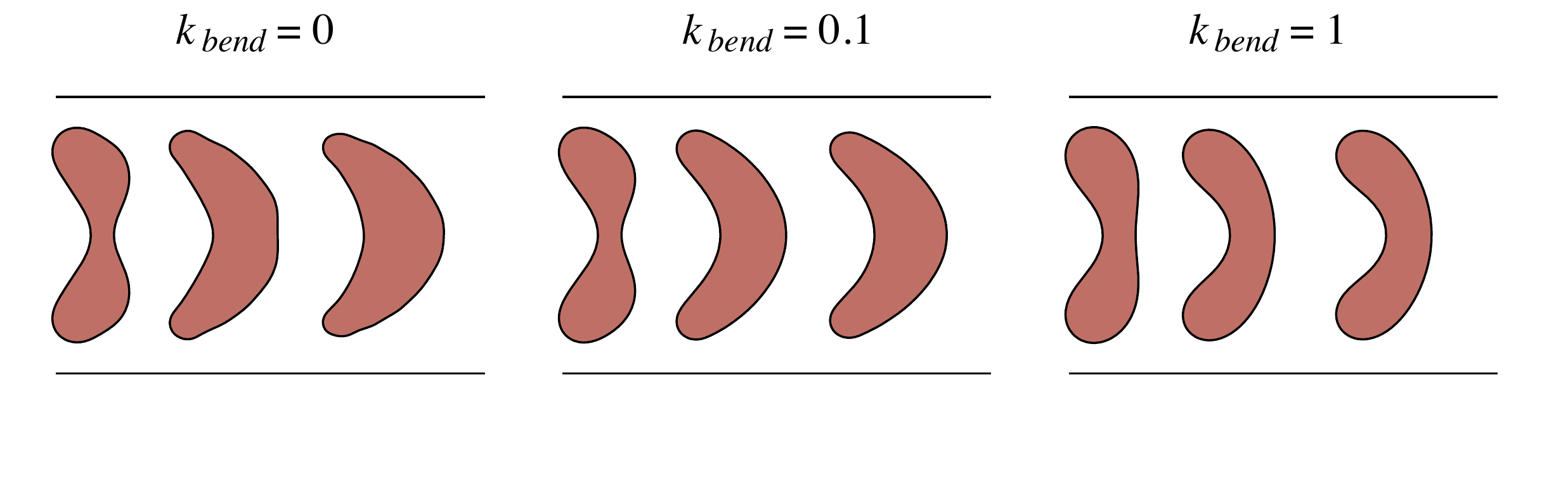}
\caption{Cross section of the cell and capillary tube through the plane $y=0$ at several time values ($t=0.005, 0.025$, and $0.05$ from left to right) for when $k_{bend}=0$ (left), 0.1 (middle), and 1 (right) pN $\cdot$ $\mu$m. The cell was simulated with SHVD ($m=400$ interpolation and $n=8281$ evaluation points).}
\label{fig:rbcbending}
\end{center}
\end{figure}

The bending rigidity is defined as the cell's resistance to being bent from a spherical shape with uniform curvature. Fig.\ \ref{fig:rbcbending} shows a slice of the RBC and capillary through the plane $y=0$ at several time points during a simulation. The cell's deformation decreases as the value of bending rigidity increases. In particular, we observe that a cell with no bending rigidity initially deforms to become pinched at its center, then flat, and then curved. Meanwhile, the front edge of the cell maintains more uniform mean curvature throughout the simulation when $k_{bend}=1$ pN $\cdot$ $\mu$m than with smaller values of $k_{bend}$. 

%
%
\subsection{Model of cellular blebbing}
\label{sec:bleb}
%
%

Blebs are round membrane protrusions that some cells use during migration \cite{CharrasPaluch}.
Blebbing occurs when a thin layer of the cytoskeleton called the cortex delaminates from the cell membrane, where it is normally attached via linker (ERM) proteins \cite{Charras2008}.
Cells that bleb are usually under high cortical tension due to myosin molecular motors pulling actin filaments with respect to each other \cite{JoannyPaluch}, a process  referred to as actomyosin contractility. 
High actomyosin contractility within the cortex leads to high intracellular pressure (on the order of $10$ Pa). 
When high intracellular pressure occurs along with a localized defect in the cortex or loss of 
membrane-cortex adhesion proteins, the membrane separates from the cortex, and cytoplasmic flow leads to bleb expansion. Hence, blebbing is a discontinuous deformation of a localized part of the cell in contrast to the global deformation of an RBC in a microcapillary. 

Several mathematical models of cellular blebbing 
have been developed to understand the interactions between the cytoskeleton, cell membrane, and adhesion proteins. These models typically make simplifying assumptions, such as the system being in steady state
\cite{JoannyPaluch,WoolleyJTB2015}, or are limited to two spatial dimensions \cite{STRYCHALSKI:BLEB:2010,StryGuyBJ2016}.
In this section, we present a three-dimensional extension of the model from \cite{STRYCHALSKI:BLEB:2010} as a first step towards understanding both the dynamics and geometry of these cellular protrusions. Our goal here is to present a dynamic 3D model of cellular blebbing where membrane forces are computed using both LDSM and SHVD. 
We compare bleb expansion dynamics and cell shapes generated by these methods. We also describe a modification to the existing method  from \cite{fai,sullivan2008curvatures} for computing mean curvature when using LDSM.
 
\subsubsection{Mathematical model}
\begin{figure}[!htbp]
\begin{center}
\includegraphics[width = 0.5\textwidth]{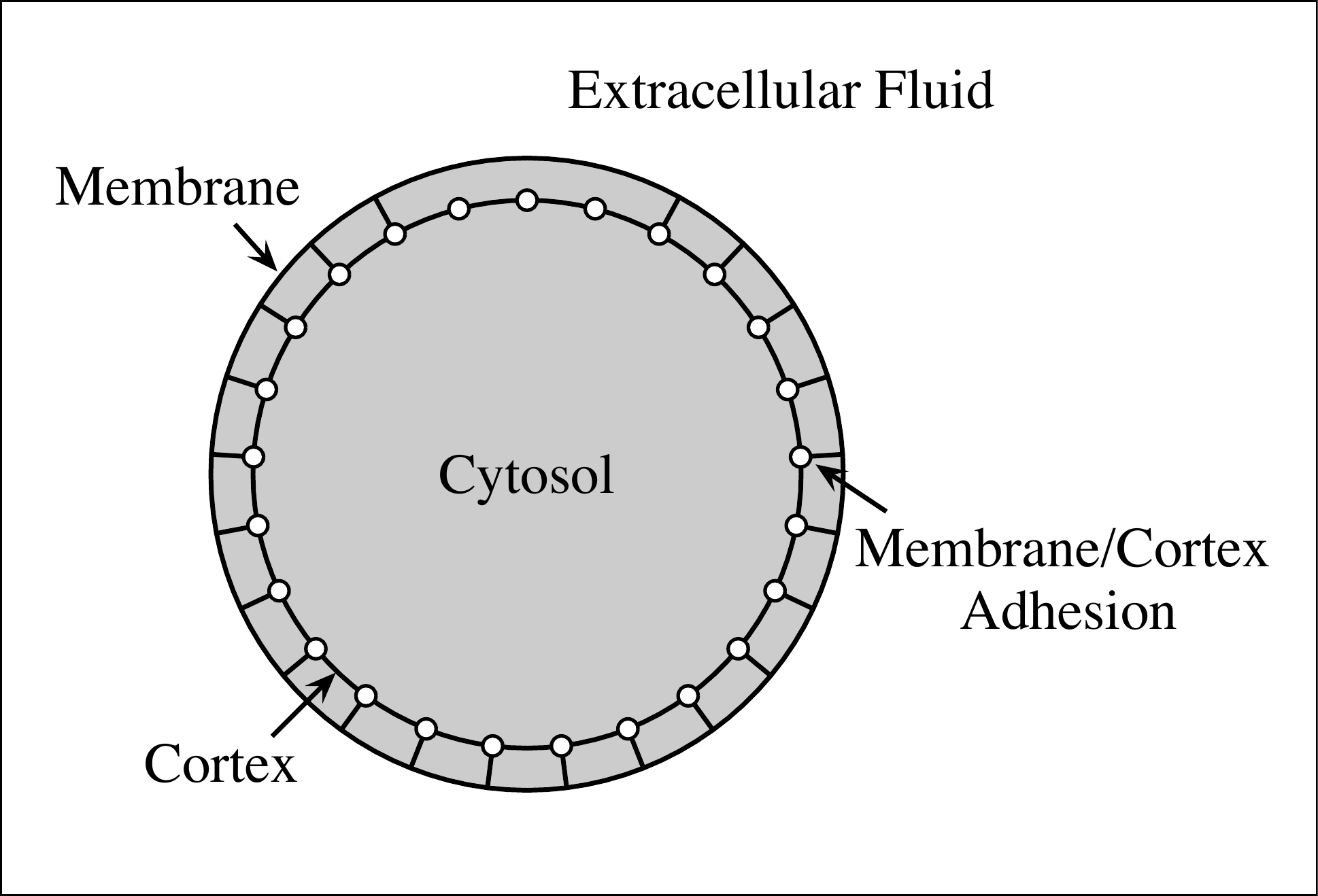}
\caption{Bleb model schematic as a slice of the cell through the $xz$-plane. A bleb is initiated by removing adhesive links in a small region at the top of the cell.}
\label{fig:blebschematic}
\end{center}
\end{figure}

Our mathematical model consists of the membrane, cortex, membrane-cortex adhesion, and intracellular fluid (see Fig.\ \ref{fig:blebschematic}). 
The intracellular fluid represents the cytosol, the liquid part of the cytoplasm. The membrane is modeled as an impermeable elastic shell that experiences forces due to surface tension and stretching. The cortex is modeled as a permeable elastic structure, and membrane-cortex adhesion is modeled by elastic springs connecting the membrane to the cortex. The cytosol is modeled as a viscous incompressible fluid governed by the fluid equations
\begin{gather}
\label{eq:blebfluid}
\mu \Delta \vec{u} - \nabla p + \vec{f}^{\rm \, mem}_{\rm elastic} + \vec{f}^{\rm \, mem/cortex}_{\rm adh} +\vec{f}^{\rm \, cortex}_{\rm drag} = \vec{0} \\
\nonumber \nabla \cdot \vec{u} = 0,
\end{gather}
where the $\vec{f}_i$'s denote spread Lagrangian force densities due to membrane elasticity, membrane-cortex adhesion, and cortical drag. The drag force on the cortex is balanced by elastic forces within the
cortex and adhesion to the membrane so that
\begin{equation}
\vec{F}_{\rm drag}^{\rm \, cortex}+ \vec{F}_{\rm elastic} ^{\rm \, cortex} +
 \vec{F}_{\rm adh}^{\rm \, cortex/mem} = \vec{0}. 
 \label{vfbleb_cortex_bal:eq}
\end{equation}
The drag force density on the fluid is related to the cortex drag
force density by
\begin{equation}
\vec{f}_{\rm drag}^{\, \rm cortex} = -\textrm{S} \vec{F}_{\, \rm drag}^{\, \rm cortex},
\end{equation}
where $\textrm{S}$ denotes the spreading operator from Eq.\ \eqref{eq:spread}.

Forces due to membrane elasticity are computed by either LDSM or SHVD with surface tension and neo-Hookean energy densities.
The cortex is modeled by a different approach, similar to the lattice-spring model described in \cite{strychalski2015poroelastic}. Specifically, the cortex is treated as a collection of elastic springs, where the force at the $i$th point is computed by 
\begin{equation}
\label{eq:3d_mem}
\vec{F}_{i} \, dA_{i}  = \sum_{j}  k^{\rm cortex}_{ij}
\Big( \frac{\norm{\vec{X}^{\rm cortex}_i-\vec{X}^{\rm cortex}_j }}{d \ell_{ij}} \Big) 
\frac{\vec{X}^{\rm cortex}_i-\vec{X}^{\rm cortex}_j}
{ \norm{\vec{X}^{\rm cortex}_i-\vec{X}^{\rm cortex}_j}},
\end{equation}
where $dA_i$ is the reference area differential, $d \ell_{ij}$ is the reference length 
of the spring connecting the $i$th point to the $j$th point on the surface, and $j$ 
is the index of all points connected to the $i$th point.
Note the area weights are the same as those from LDSM in Eq.\ \eqref{eq:faiwts}.
The spring coefficient connecting the $i$th point to the $j$th point is calculated as 
\begin{equation}
k^{\rm cortex}_{ij} = \frac{8k_{ij}}{3 d \ell_{ij}} \left(\frac{dA_i+dA_j}{2} \right),
\end{equation}
as in the lattice-spring model from \cite{strychalski2015poroelastic}. This particular model 
of the cortex was used so that cortical ablation experiments could 
easily be performed by setting $k^{\rm cortex}_{ij}=0$ in a localized region. We plan to explore this mechanism of bleb initiation in a future publication.

Membrane-cortex adhesion is modeled by elastic springs attaching the membrane to the
cortex with a force density given by
\begin{equation}
\label{eq:mem_cortex_adh}
\vec{F}_{\rm adh}^{\, \rm mem/cortex} = k_{\rm adh}^{\rm mem/cortex} 
\Big( \norm{\vec{X}_{\rm mem}-\vec{X}_{\rm cortex}} \Big) 
\frac{ \vec{X}_{\rm mem}-\vec{X}_{\rm cortex}}
{ \norm{\vec{X}_{\rm mem}-\vec{X}_{\rm cortex}}}.
\end{equation}
The membrane-cortex adhesion stiffness coefficient $k_{\rm adh}^{\rm cortex/mem}$ 
was chosen so that the velocity of the cortex is zero if no bleb is initiated.
The adhesion force density on the membrane 
is the opposite of the corresponding force density on the cortex, with the proper scaling to 
ensure that the two forces balance,
\begin{equation}
\int \limits_{\Omega} \textrm{S} \vec{F}_{\rm adhesion}^{\, \rm mem/cortex} \, d \vec{x} + 
\int \limits_{\Omega} \textrm{S} \vec{F}_{\rm adhesion}^{\, \rm cortex/mem} \, d \vec{x}  = 0.
\end{equation}
Given the stiffness coefficient $k_{\rm adh}^{\rm cortex/mem}$, 
 the corresponding stiffness coefficient for 
the cortex is obtained by $k_{\rm adh}^{\rm mem/cortex} = 
k_{\rm adh}^{\rm cortex/mem} dA^{\rm cortex}_i/dA^{\rm mem}_i$, where 
$dA_i$ represents the area differential in reference coordinates. For the membrane, these $dA_i$ values are given by Eq.\ \eqref{eq:wtseq} for SHVD and Eq.\ \eqref{eq:faiwts} for LDSM. 

Given a configuration of the membrane and cortex, the Lagrangian force densities are computed, and the velocities of the fluid and cortex are obtained by
solving
Eqs.\ \eqref{eq:blebfluid} and \eqref{vfbleb_cortex_bal:eq}. 
The positions of the
membrane and cortex are then updated with their respective velocities,
\begin{align}
\frac{d \vec{X}_{\rm mem}}{dt} &= \textrm{S}^*\vec{u} = \vec{U}, \\
\label{Eq:CortexVelocity} \frac{d \vec{X}_{\rm cortex}}{dt}  &  = 
\frac{1}{\xi} \lf(\vec{F}_{\rm elastic }^{\, \rm cortex}+ \vec{F}_{\rm attach}^{\, \rm cortex/mem} \rt) + 
\vec{U} =  \vec{U}_{\rm cortex},
\end{align}
where $\xi$ is the drag coefficient of the cortex.
A bleb is initiated by removing adhesive links between the membrane and cortex in a small region given by the equation $\phi \leq 2\pi/25$ in spherical coordinates (see Fig. \ref{fig:blebschematic}).
Model parameters values are provided in 
Table \ref{ModelParametersNoCyto}.
%
%
\begin{table}
\centering
\begin{tabular}{|lllc|} 
\hline
Symbol & Quantity & Value & Source \\
\hline
$r_{\rm mem}$ & Cell radius & 10 $\mu$m & \cite{JoannyPaluch}\\ 
$r_{\rm cortex}$ & Cortex radius & 9.99 $\mu$m & \cite{StryGuyBJ2016} \\
$\gamma_{\rm mem}$ & Membrane surface tension& 40 pN/$\mu$m & \cite{JoannyPaluch}\\
$\kappa_{E}$ & Membrane bulk modulus & 40 pN/$\mu$m & \cite{StryGuyBJ2016}  \\
$\mu_{E}$ & Membrane shear modulus & 40 pN/$\mu$m & \cite{StryGuyBJ2016}  \\
$k_{\rm cortex}$ & Cortical  stiffness coefficient (average value) & 87 pN/$\mu$m$^2$ &   \\
$k_{\rm adh}^{\rm mem/cortex}$ & Membrane/cortex adhesion & $3 \cdot 10^4$ pN/$\mu$m$^3$ &   \\
& stiffness coefficient &  & \\
$\mu$ & Cytosolic viscosity & $1$ Pa\--s &  \cite{Charras2008,JoannyPaluch,Moeendarbary2013} \\
$\xi$ & Cortical drag coefficient & 10 pN-s/$\mu$m$^3$  & \cite{STRYCHALSKI:BLEB:2010} \\ 
\hline
\end{tabular}
\caption{Parameters for the blebbing model.}
\label{ModelParametersNoCyto}
\end{table}

The model is simulated on a box with volume 
$30^3 $ $\mu \textrm{m}^3$ and 
$\Delta x=\Delta y=\Delta z = 30/32.$ The membrane and cortex are initialized to be spheres centered at the point $(15,15,15)$ with 4761 Lagrangian points, corresponding to approximately 2 Lagrangian points per Eulerian grid cube. For the SHVD model, we compute all solutions with 4761 evaluation points and vary the number of interpolation points to be 625, 1024, and 4761. The lattice spring model on the cortex results in large spurious velocities at the beginning of the simulation. 
Before initiating a bleb, we allow the forces in the model to equilibrate, meaning the velocity on the boundary decreases to approximately 0.01 $\mu$m/s. This process corresponds to 0.2 s of simulation time. Results are reported in time units after equilibration.
\subsubsection{Model results}
\begin{figure}[htbp]
\begin{center}
\includegraphics[width = 0.9\textwidth]{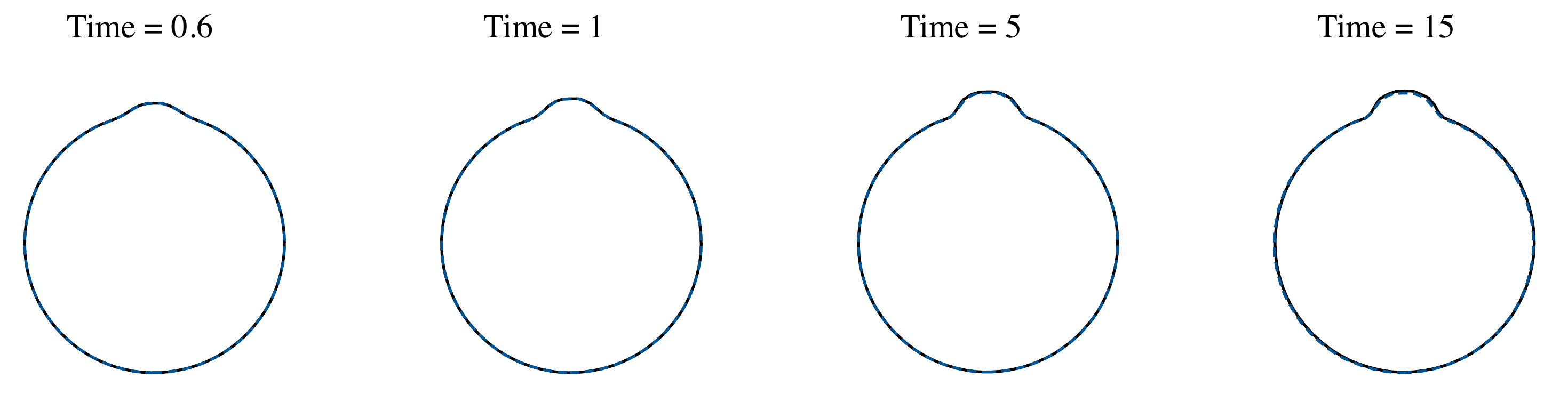}
\caption{Slices through the cell membrane by the $yz$ plane when $x=15$ at several time values. The black line denotes the membrane position computed by LDSM, and the blue dashed line indicates the membrane position computed by SHVD with 1024 interpolation points and 4761 evaluation points.}
\label{fig:blebslices}
\end{center}
\end{figure}

Pressure is initially uniform inside the cell. When adhesive links between the membrane and cortex are removed, cortical tension is no longer transmitted to the membrane locally, resulting in a pressure gradient and bleb expansion. Fig.\ \ref{fig:blebslices} shows cross sections of the cell membrane as a bleb expands using both LSDM and SHVD with 1024 interpolation points. The cross sections show that the models have good agreement with each other. 

\begin{figure}[htbp]
\begin{center}
\includegraphics[width = 0.9\textwidth]{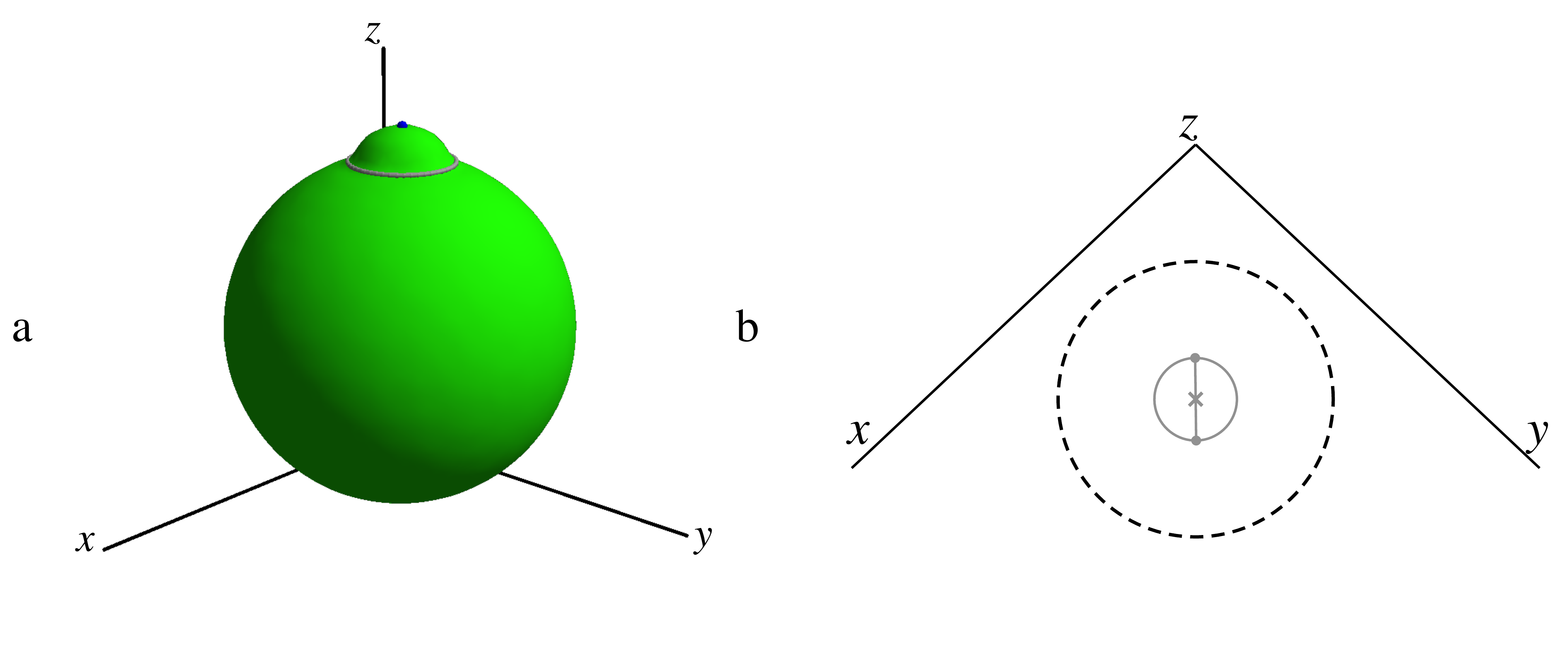}
\caption{The bleb ring is defined by the boundary between the region where adhesive links connect the membrane to the cortex and the region where the links are removed, indicated by the gray ring in (a). The membrane shell is graphed in green in (a). Bleb diameter is found by finding the point on the bleb ring that is closest to the point (0,0,30) and the point directly across from it. A top town view of the cell and bleb ring is shown in (b) with the two points that define the bleb's diameter. The dashed black line indicates the membrane and the solid gray circle indicates the bleb ring. The distance from the midpoint of the bleb diameter (indicated by $\times$ in the top-down view) to the point with the largest z-coordinate value, as indicated by the blue point in (a), minus its initial value defines bleb size.}
\label{fig:blebmidpointsize}
\end{center}
\end{figure}

In order to obtain a localized measurement of bleb size, we compute the bleb ring as the one dimensional boundary between the region where the membrane is attached to the cortex and the region where the membrane is detached from it (see Fig.\ \ref{fig:blebmidpointsize}). Bleb diameter is measured as the distance from the point closest to the $z$-axis to the point farthest away from it on the bleb ring. The midpoint of this line is computed, and bleb size is defined to be the vertical distance from the midpoint of the bleb diameter line to the point on the bleb with the largest $z$-coordinate (blue point in Fig.\ \ref{fig:blebmidpointsize}(a)) after subtracting the initial value. Note that this is a 3D extension of the method described in \cite{StryGuyBJ2016}.

\begin{figure}[htbp]
\begin{center}
\includegraphics[width = 0.67\textwidth]{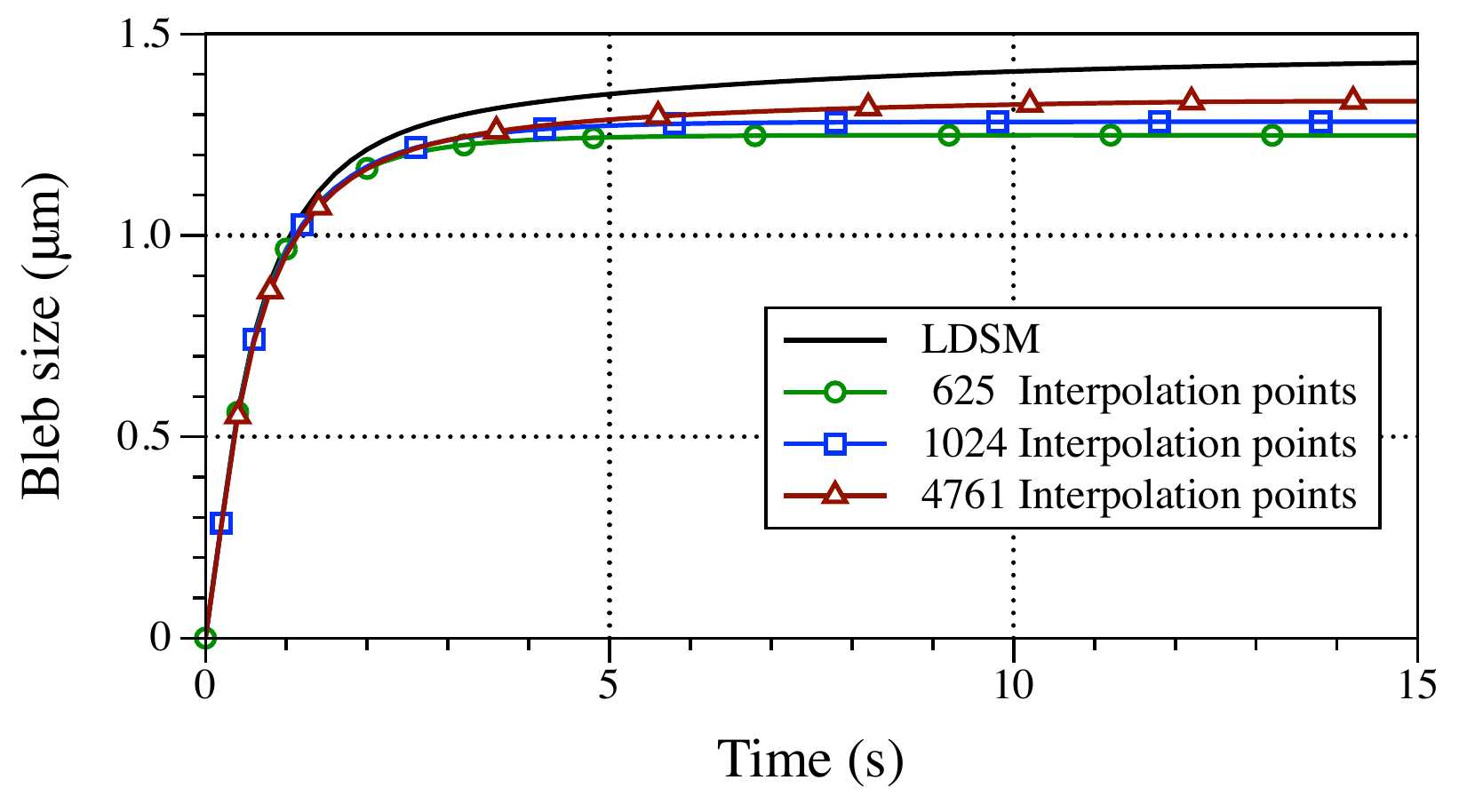}
\caption{Bleb size over time when the membrane position is computed by LDSM (black line) and by SHVD with various numbers of interpolation points and 4761 evaluation points.}
\label{fig:blebsizetime}
\end{center}
\end{figure}

Bleb size over time when membrane forces are computed by LDSM and with several different values of interpolation points in SHVD is shown in Fig.\ \ref{fig:blebsizetime}.
The data show very similar bleb expansion dynamics and steady state bleb size for all simulations.
We observe that the bleb size increases as the bleb neck is better resolved with a larger number of interpolation points. Once we reach $m=n=4761$ (red triangles in Fig.\ \ref{fig:blebsizetime}), each evaluation is updated explicitly with the local fluid velocity. In this case (when $m=n$), the difference between our method and LDSM is due only to the method of computing forces.

\begin{figure}[htbp]
\begin{center}
\includegraphics[width = 1\textwidth]{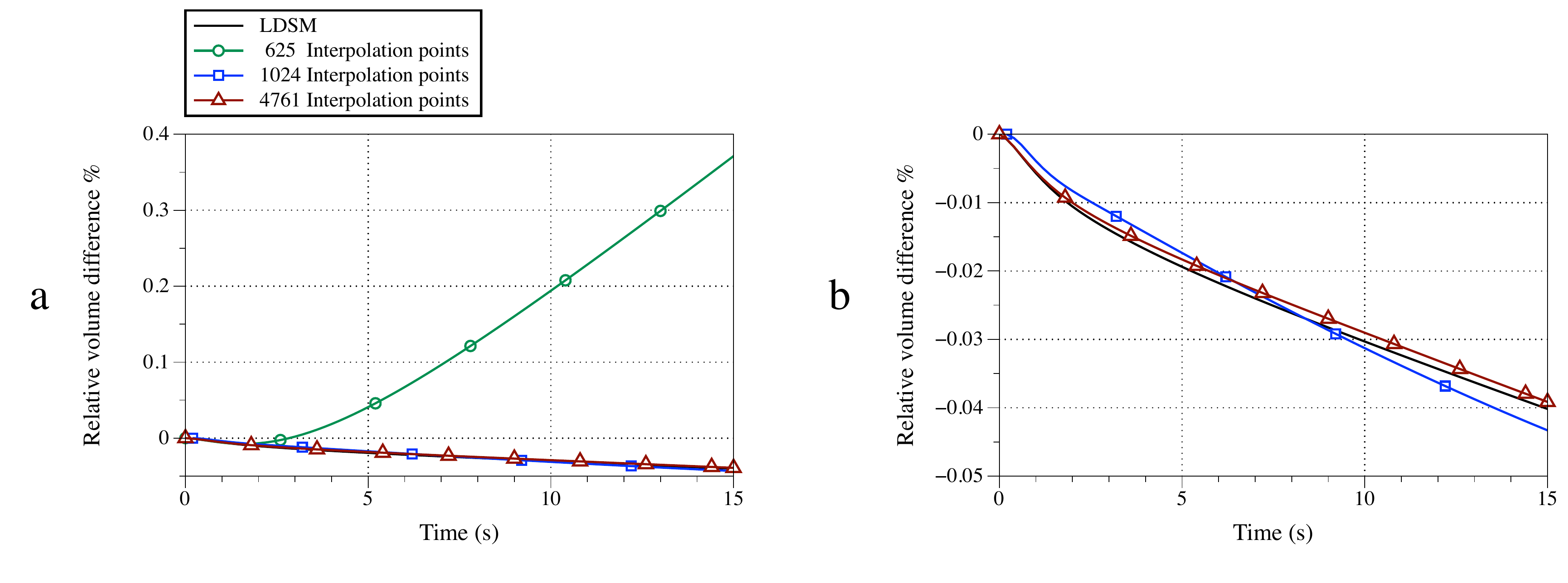}
\caption{Membrane volume over time when SHVD is used for the membrane with 625, 1024, and 4761 interpolation points. Data from the membrane modeled with LSDM are shown in black for comparison. The graph on the right is a zoomed in view of the same data from the left.}
\label{fig:blebvol}
\end{center}
\end{figure}

Bleb volume over time is shown in Fig.\ \ref{fig:blebvol}, 
where percentage relative volume difference is defined as
$(V(t)-V(0))/V(0) \times 100$, where $V(t)$ is the volume of the membrane at time $t$. Results from LDSM and SHVD with $1024$ and $4761$ interpolation points show good agreement. Volume actually increases over time with the SHVD model using 625 interpolation points. Because of the sharp transition from 
a region with membrane-cortex adhesion to one without adhesion, a large number of interpolation points is needed to resolve the region near the bleb neck. 

Notably, the number of overall interpolation points required to resolve the bleb neck, which is a highly localized region of deformation, is higher than the number of points required to resolve the deformed RBC in Section \ref{sec:RBCs} because the points overall must be evenly spread over the cell to yield an accurate interpolant. Hence more points must be added throughout the cell, not just in the localized region, in order to increase resolution of the bleb neck. Unlike our RBC application, adding points up to the limit $m=n$ actually \textit{stabilizes} the steady state membrane shapes. We believe this is because most of the membrane remains still while tethered to the cortex. Both of these factors (the lack of motion from a background flow and the tether forces on the membrane) prevent higher order modes from distorting the membrane shape.

\begin{figure}[!tbp]
\begin{center}
\includegraphics[width = 0.7\textwidth]{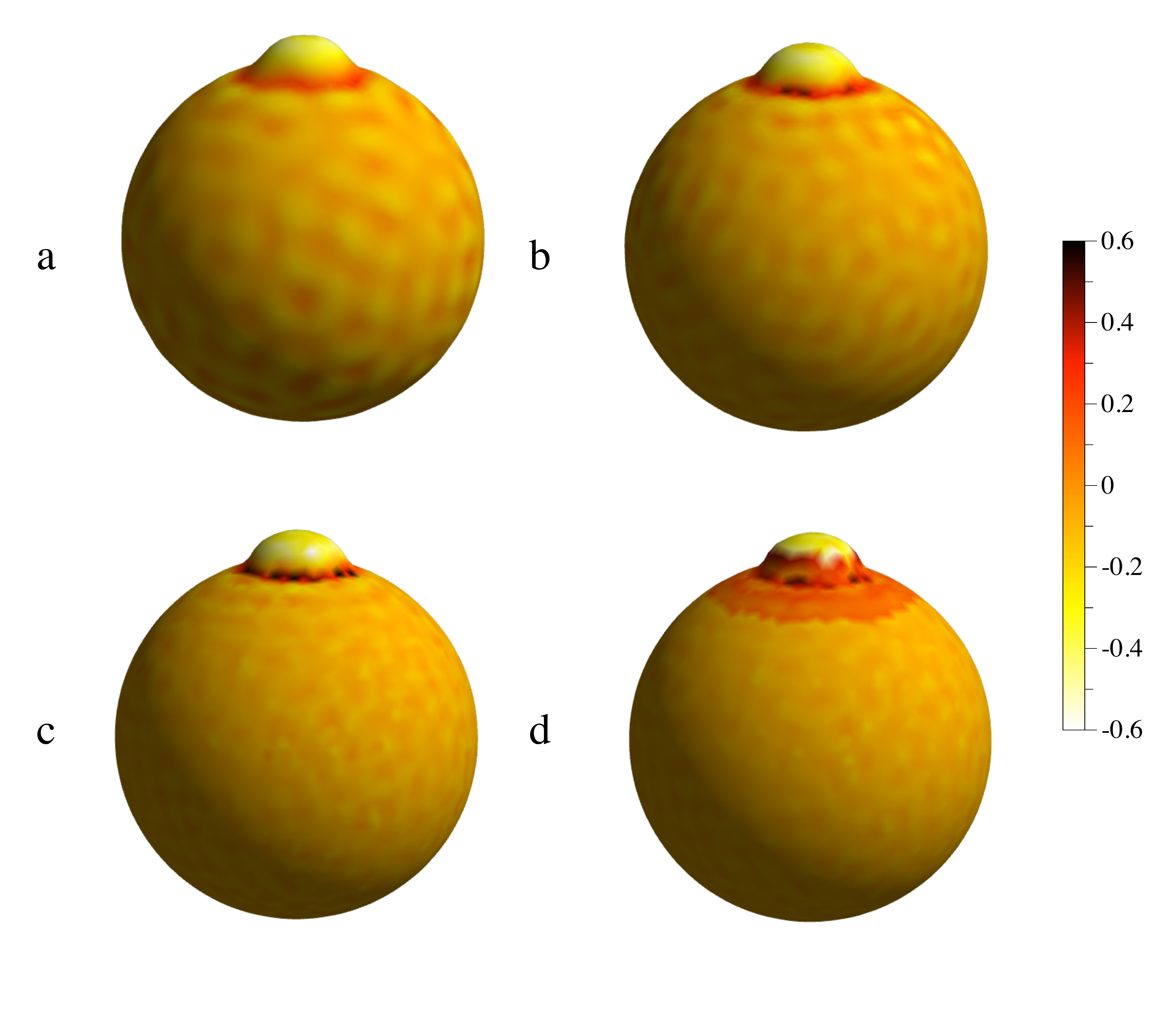}
\caption{Membrane mean curvature after bleb expansion at $t=10$ s when the membrane is modeled using SHVD with (a) 625, (b) 1024, and (c) 4761 interpolation points. Mean curvature from LSDM with 4761 points is shown in (d).}
\label{fig:blebcurveSH}
\end{center}
\end{figure}

Membrane mean curvature after bleb expansion (at $t=15$ s) is shown in Fig.\ \ref{fig:blebcurveSH} for the same data from Figs.\ \ref{fig:blebsizetime} and \ref{fig:blebvol}.  In order to compute curvature for LDSM, we follow the method described in \cite{fai,sullivan2008curvatures}, where curvature is computed on edges and vertices of surface triangles. This method results in curvatures that are always positive. In order to obtain signed mean curvature, we post-process the surface to compute its convex hull. Curvature values at points on the convex part of the surface are assigned to be negative while the rest of the curvature values (on the concave part of the surface) maintain their positive sign.
In SHVD with 625 interpolation points, the bleb appears more broad which indicates it is not well-resolved. The membrane is effectively smoothed out if not enough interpolation points are used, and volume slightly increases.
The membrane appears more dimpled in SHVD with 625 interpolation points compared to data from the other simulations. This can be partially attributed to the adhesion model, where membrane evaluation points are connected to points on the cortex. Since membrane evaluation points are not explicitly updated with the fluid velocity, adhesion forces are slightly perturbed when compared to model results from Fig.\ \ref{fig:blebcurveSH} (c) and (d), where the interpolation and evaluation points are updated with the local fluid velocity (although some dimpling is present even in Fig.\ \ref{fig:blebcurveSH} (c) and (d) due to membrane-cortex adhesion). The bleb appears to be well-resolved in SHVD with 1024 interpolation points with a smooth membrane.  Membrane curvature from LSDM appears discontinuous and noisy when compared to simulation data from SHVD. We emphasize that SHVD more accurately computes curvature than LDSM given the ability of the spherical harmonic interpolant to capture the membrane shape. 

The results from this section show that cell shape and bleb expansion dynamics are qualitatively similar when the membrane is modeled by LDSM and SHVD. If accurate geometric information, such as membrane curvature, is desired from a simulation, SHVD should be used for computing membrane forces. In \cite{STRYCHALSKI:BLEB:2010}, the authors used a 2D dynamic model to determine the role of cortical drag and cytoplasmic viscosity on bleb expansion dynamics. We find that for obtaining dynamic information such as bleb expansion time, LDSM is sufficiently accurate.
In the next section, we comment on the computational cost associated with SHVD when compared to LDSM.
%

%
%
\section{Discussion}
\label{sec:discuss}
%
%

We show that by representing surfaces with spherical harmonic interpolants, force densities can be computed more accurately as continuous functions by taking the variational derivative of an energy density functional than by discretizing the surface first into piecewise planar elements. Furthermore, geometric information such as curvature is more accurately computed with the continuous representation. We presented tests and applications for zero Reynolds number flow, but it is straightforward to extend our method to fluid-structure interaction problems where the fluid is governed by the Navier-Stokes equations. This is an advantage of our method compared to boundary integral methods that employ variational methods for computing forces, in that ours can also be used when modeling red blood cells in intermediate Reynolds number flows, for example \cite{fai}.

A limitation of our method, as with all IB methods, is its overall first order accuracy in space because the regularized delta function used in the spreading operator does not accurately capture discontinuities in the pressure and derivatives of the velocity across the infinitely thin interfaces. Immersed interface methods \cite{li2006immersed} have been proposed to improve the accuracy of the IB method. These ``IM methods'' depend on quantifying the jump conditions at the interface, which are functions of geometric information, Lagrangian forces, and their derivatives \cite{xu20083d}. The order of the method depends in part on how accurately these two quantities are computed. Because our method is spectrally accurate in both, a future goal of ours is to implement it in an IM framework, thereby increasing the overall accuracy of the method.

We also provide a method that ensures the integral of the Lagrangian force density is equal to the integral of the spread force density. 
This allows us to apply our method to models of structures immersed in Stokes flow, a regime that is prevalent across most cellular-level applications. We do note that our formulation relies on the structure being a closed surface.

Our method is applied to models of red blood cells in capillary flow and cellular blebbing. We have intentionally chosen two applications with very different dynamics and behavior in order to showcase the versatility of our method and highlight its strengths and limitations. In the RBC model, we see a smooth change in the entire surface over time. This allows us to resolve the shape with $\mathcal{O}(100)$ interpolation points. Lower numbers of interpolation points do not resolve the shape well enough, and higher numbers lead to numerical instabilities as higher frequency modes amplify non-physical bumps and noise on the cell surface. 

Meanwhile, cellular blebbing involves a discontinuous change in a small region of the cell due to a discontinuity in the membrane-cortex adhesion force density. In our present formulation, the points must be evenly distributed to build a reliable interpolant and a large number of points must be used on the entire cell to accurately resolve the bleb neck. Therefore it is necessary to choose $m$ (number of interpolation points) and $n$ (number of evaluation points) to be on the same order of magnitude. We found that measuring volume of the membrane over time is an indicator of surface resolution; if membrane volume increases during a simulation, the surface is likely under-resolved and more interpolation points are needed.
We hypothesize that the physical linkage of the membrane to the cortex stabilizes the higher frequency modes in the spherical harmonic interpolant, and numerical instabilities from a large number of interpolation points are not observed in this model.

The restriction $m \ll n$ is for computational efficiency; our method is most efficient if the $\mathcal{O}(mn)$ cost of computing the forces is comparable to the $\mathcal{O}(n)$ cost in LDSM. This is the case for the RBC model, where the surface is well-resolved with $\mathcal{O}(100)$ interpolation points and the computational time of our method is of the same order of magnitude as LDSM (approximately 1 minute). However, to simulate the blebbing model, $m$ must be much larger, and the computational time to compute the forces increases significantly. 
A simulation of bleb expansion took approximately $10$ hours with SHVD compared to approximately 1 hour for LDSM. While our model is versatile enough to be applied across a wide variety of applications, this additional cost should be taken into account.

We plan to use these principles to guide our future work. For example, our method would be well-suited for simulating a model of cytokinesis, the division of the cytoplasm when one cell divides into two cells \cite{AlbertsCell}. 
During cytokinesis, the entire membrane deforms, and our method allows the force density to be computed at any point on the surface. Models of microtubule attachment to the cortex during spindle positioning will benefit from knowing the elastic force exactly at the points where the microtubules attach to the cortex. The noisy force data from LDSM could lead to inaccurate results. Finally, some models assume that spontaneous membrane curvature initializes furrow ingression during cytokinesis \cite{Schlomovitz08}. Our framework provides an accurate method for computing curvature in a dynamic 3D model of the early stages of cytokinesis.

Other future work includes extending our framework to include adaptive meshing by adding interpolation points when the surface shape undergoes large local deformations and surface features become unresolved. This capability, along with finding a more efficient way to compute the force densities (e.g. in parallel), could make our method more suitable for use with the blebbing model and allow us to investigate different bleb initiation mechanisms using our method. The formulation we have outlined can be applied to other interpolating basis functions, including radial basis functions. Because spherical harmonic basis functions were able to provide exact force calculations for very smooth (spherical) shapes, our method was computationally advantageous when it was applied to shapes  resembling spheres. We anticipate that accuracy for arbitrary shapes (including blebbing cells) could be achieved by shifting the basis functions according to the application. We plan to investigate these topics in future publications.

\section*{Acknowledgements}
This work was supported in part by grant \#429808 from the Simons Foundation to W.S. A.K. was supported by a National Science Foundation grant DMS-1521748. The authors would like to thank Grady B. Wright for discussions about quadrature weights on the sphere, Varun Shankar for helpful discussions about surface representation, and Boyce Griffith for suggestions about computing forces due to bending. The authors also thank Thomas G. Fai and Charles S. Peskin for notes about modeling neo-Hookean materials in curvilinear coordinates.
\appendix
\section{Derivation of Force Density, $\vec{F}$}
%
%
Here we derive in detail the force density in Eq.\ \eqref{eq:vdforce} from the variational derivative of the energy density functional. We make use of dyad expansions of tensors ($\vec{e_i}\vec{e_j}$ denotes the $i$th row and $j$th column of the given tensor) and the definitions of the dot and double dot products ($\vec{e_i} \cdot \vec{e_j} = 1$ if $i=j$ and 0 otherwise, $\vec{e_i}\vec{e_j} \cddot \vec{e_k}\vec{e_p} = 1$ if $j=k$ and $i=p$ and 0 otherwise). 

\subsection{Dyad expansion of $\ten{P}$ allows us to recover Eq.\ \eqref{eq:frechsum} from Eq.\ \eqref{eq:deldotP}}
We begin with showing the equivalence of Eq.\ \eqref{eq:frechsum} and Eq.\ \eqref{eq:deldotP}. In Eq.\ \eqref{eq:frechsum}, the variational derivative is given by
\begin{equation}
\frac{\delta \mathcal{E}} {\delta \vec{X}} = 
-\sum_{i=1}^3 \sum_{j=1}^2 {\frac{1}{\sqrt{\det{\Gz}}}\frac{\partial (\ten{P}_{ij} \sqrt{\det{\Gz}})}{\partial q_j}}.
\end{equation}
Considering the product in Eq.\ \eqref{eq:deldotP},
\begin{align}
\nonumber -\nabla_{\theta} \cdot \ten{P}^{\textrm{T}} & = -\frac{1}{\sqrt{\Det{\G_0}}} \left(\frac{\partial}{\partial \lambda} \sqrt{\Det{\G_0}} \textrm{ } \vec{e_1} + \frac{\partial}{\partial \theta} \sqrt{\Det{\G_0}} \textrm{ } \vec{e_2}\right)\cdot \\[4 pt]  & (\ten{P}_{11}\vec{e_1}\vec{e_1} + \ten{P}_{21}\vec{e_1}\vec{e_2} +
\ten{P}_{31}\vec{e_1}\vec{e_3} + \ten{P}_{12}\vec{e_2}\vec{e_1} + \ten{P}_{22}\vec{e_2}\vec{e_2} + \ten{P}_{32}\vec{e_2}\vec{e_3})\\[4 pt]
& = -\sum_{i=1}^3 {\frac{1}{\sqrt{\det{\Gz}}}\left(\frac{\partial}{\partial \lambda}\left(\ten{P}_{i1} \sqrt{\det{\Gz}}\right)+\frac{\partial}{\partial \theta}\left(\ten{P}_{i2} \sqrt{\det{\Gz}}\right)\right)}\\[4 pt]
\label{eq:Alastexpansion}
& = -\sum_{i=1}^3 \sum_{j=1}^2{\frac{1}{\sqrt{\det{\Gz}}}\left(\frac{\partial}{\partial q_j}\left(\ten{P}_{ij} \sqrt{\det{\Gz}}\right)\right)}.
\end{align}
Eq.\ \eqref{eq:Alastexpansion} is the same expression as that derived from the variational derivative definition in Eq.\ \eqref{eq:frechsum},
\begin{equation}
  \vec{F}=-\nabla_{\theta} \cdot \ten{P}^{\textrm{T}}.
\end{equation}
\subsection{Chain rule expansion of $\ten{P}^{\textrm{T}}$}
We here prove Eq.\ \eqref{eq:tenchain}, 
\begin{equation}
\ten{P}^{\textrm{T}}=\frac{\partial{W}}{\partial{(\mat{\nabla_{\theta}X})^\textrm{T}}} = \frac{\partial W}{\partial \gamma} \cddot \frac{\partial \gamma}{\partial (\mat{\nabla_{\theta}X})^\textrm{T}} = \ten{S} \cddot \frac{\partial \gamma}{\partial (\mat{\nabla_{\theta}X})^\textrm{T}} = \ten{S} \cdot (\mat{\nabla_{\theta}X})^\textrm{T}.
\end{equation}
We begin with the double dot product $\fraco{\partial W}{\partial \gamma} \cddot \fraco{\partial \gamma}{\partial (\mat{\nabla_{\theta}X})^\textrm{T}}$. Recall from Eq.\ \eqref{eq:GLtensor} that $\gamma$ is the Green-Lagrange strain energy tensor, given by $\fraco{1}{2}(\G-\Gz)$. Because $\gamma$ is symmetric, the tensor $\fraco{\partial W}{\partial \gamma} = \ten{S}$ (the second Piola-Kirchhoff stress tensor) is also symmetric, and its dyad expansion is
 \begin{equation}
 \frac{\partial W}{\partial \gamma} = \ten{S} = \ten{S}_{11} \vec{e_1}\vec{e_1} + \ten{S}_{12}\vec{e_1}\vec{e_2} + \ten{S}_{21}\vec{e_2}\vec{e_1} + \ten{S}_{22}\vec{e_2}\vec{e_2}.
 \end{equation}
  Now, consider the fourth order tensor $\fraco{\partial \gamma}{\partial (\mat{\nabla_{\theta}X})^\textrm{T}}$. We denote 
 \begin{equation*}
 \left(\frac{\partial \gamma}{\partial (\mathbf{\nabla_{\theta}X})^\textrm{T}}\right)_{ijk\ell} = \frac{\partial \gamma_{ij}}{\partial (\mathbf{\nabla_{\theta}X})^\textrm{T}_{k\ell}}.
 \end{equation*}
 We write $\gamma$ as
 \begin{equation}
 \gamma  =  \frac{1}{2}\begin{pmatrix}
\fraco{\partial\bm{X}}{\partial \lambda} \cdot \fraco{\partial\bm{X}}{\partial \lambda} - \fraco{\partial\bm{Z}}{\partial \lambda} \cdot \fraco{\partial\bm{Z}}{\partial \lambda} & \fraco{\partial\bm{X}}{\partial \lambda} \cdot \fraco{\partial\bm{X}}{\partial \theta} - \fraco{\partial\bm{Z}}{\partial \lambda} \cdot \fraco{\partial\bm{Z}}{\partial \theta}\\[10 pt]
\fraco{\partial\bm{X}}{\partial \theta} \cdot \fraco{\partial\bm{X}}{\partial \lambda} - \fraco{\partial\bm{Z}}{\partial \theta} \cdot \fraco{\partial\bm{Z}}{\partial \lambda} & \fraco{\partial\bm{X}}{\partial \theta} \cdot \fraco{\partial\bm{X}}{\partial \theta} - \fraco{\partial\bm{Z}}{\partial \theta} \cdot \fraco{\partial\bm{Z}}{\partial \theta}\\
        \end{pmatrix}.
\end{equation}
Expanding entrywise,
\begin{equation}
\label{eq:entrywisegam}
        \gamma = \frac{1}{2}\begin{pmatrix}
\fraco{\partial X_1^2}{\partial \lambda} + \fraco{\partial X_2^2}{\partial \lambda} + \fraco{\partial X_3^2}{\partial \lambda} - \dots  & \fraco{\partial X_1}{\partial \lambda} \fraco{\partial X_1}{\partial \theta} + \fraco{\partial X_2}{\partial \lambda} \fraco{\partial X_2}{\partial \theta} + \fraco{\partial X_3}{\partial \lambda} \fraco{\partial X_3}{\partial \theta} - \dots\\[10 pt]
\fraco{\partial X_1}{\partial \lambda} \fraco{\partial X_1}{\partial \theta} + \fraco{\partial X_2}{\partial \lambda} \fraco{\partial X_2}{\partial \theta} + \fraco{\partial X_3}{\partial \lambda} \fraco{\partial X_3}{\partial \theta} - \dots & \fraco{\partial X_1^2}{\partial \theta} + \fraco{\partial X_2^2}{\partial \theta} + \fraco{\partial X_3^2}{\partial \theta} - \dots
        \end{pmatrix}.
 \end{equation}
Eq.\ \eqref{eq:entrywisegam} allows us to easily compute derivatives with respect to the entries of $(\mathbf{\nabla_{\theta}X})^\textrm{T}$. For example, 
 \begin{align}
 \left(\frac{\partial \gamma}{\partial (\mathbf{\nabla_{\theta}X})^\textrm{T}}\right)_{ij1\ell} & = \fraco{\partial \gamma}{\partial \left(\fraco{\partial X_\ell}{\partial \lambda}\right)} = \begin{pmatrix}
 \fraco{\partial X_\ell}{\partial \lambda} & \fraco{1}{2}\fraco{\partial X_\ell}{\partial \theta}\\[6 pt]
 \fraco{1}{2}\fraco{\partial X_\ell}{\partial \theta} & 0
 \end{pmatrix} \hspace{0.2cm} 
 \end{align}
 and
 \begin{align}
 \left(\fraco{\partial \gamma}{\partial (\mathbf{\nabla_{\theta}X})^\textrm{T}}\right)_{ij2\ell} & = \fraco{\partial \gamma}{\partial \left(\fraco{\partial X_\ell}{\partial \theta}\right)} = \left(\begin{matrix}
 0 & \fraco{1}{2}\fraco{\partial X_\ell}{\partial \lambda}\\[10 pt]
 \fraco{1}{2}\fraco{\partial X_\ell}{\partial \lambda} & \fraco{\partial X_\ell}{\partial \theta}
 \end{matrix}\right).
 \end{align} 
The full expansion is
\begin{align}
\label{eq:fullexpansion}
   \frac{\partial W}{\partial \gamma} \cddot \frac{\partial \gamma}{\partial (\mathbf{\nabla_{\theta}X})^\textrm{T}} = \ten{S} \cddot \frac{\partial \gamma}{\partial (\mathbf{\nabla_{\theta}X})^\textrm{T}} = \sum_{i,j} \ten{S}_{ij} \vec{e_i} \vec{e_j} \cddot \sum_{m,p,k,\ell}  \left(\frac{\partial \gamma}{\partial (\mathbf{\nabla_{\theta}X})^\textrm{T}}\right)_{mpk\ell} \vec{e_m} \vec{e_p} \vec{e_k} \vec{e_\ell} .
\end{align}
Recall that $\vec{e_i} \vec{e_j} \cddot \vec{e_m} \vec{e_p} = 1$ if $j=m$ and $i=p$ and is 0 otherwise. This reduces the $(k,\ell)$ entry of Eq.\ \eqref{eq:fullexpansion} to
\begin{equation}
\left(\ten{S} \cddot \frac{\partial \gamma}{\partial (\mathbf{\nabla_{\theta}X})^\textrm{T}}\right)_{k\ell} =\sum_{i,j}  \ten{S}_{ij}\left(\frac{\partial \gamma}{\partial (\mathbf{\nabla_{\theta}X})^\textrm{T}}\right)_{jik\ell},
\end{equation}
or in matrix form, 
\begin{align}
& \ten{S} \cddot \frac{\partial \gamma}{\partial (\mathbf{\nabla_{\theta}X})^\textrm{T}}=\\\nonumber
& \begin{pmatrix}
\ten{S}_{11}\fraco{\partial X_1}{\partial \lambda} + \fraco{1}{2}\ten{S}_{12}\fraco{\partial X_1}{\partial \theta} + \fraco{1}{2}\ten{S}_{21}\fraco{\partial X_1}{\partial \theta} & \ten{S}_{11}\fraco{\partial X_2}{\partial \lambda} + \fraco{1}{2}\ten{S}_{12}\fraco{\partial X_2}{\partial \theta} + \fraco{1}{2}\ten{S}_{21}\fraco{\partial X_2}{\partial \theta} & \ten{S}_{11}\fraco{\partial X_3}{\partial \lambda} + \fraco{1}{2}\ten{S}_{12}\fraco{\partial X_3}{\partial \theta} + \fraco{1}{2}\ten{S}_{21}\fraco{\partial X_3}{\partial \theta} \\[10 pt]
\fraco{1}{2}\ten{S}_{12}\fraco{\partial X_1}{\partial \lambda} + \fraco{1}{2}\ten{S}_{21}\fraco{\partial X_1}{\partial \lambda} + \ten{S}_{22}\fraco{\partial X_1}{\partial \theta} & 
\fraco{1}{2}\ten{S}_{12}\fraco{\partial X_2}{\partial \lambda} + \fraco{1}{2}\ten{S}_{21}\fraco{\partial X_2}{\partial \lambda} + \ten{S}_{22}\fraco{\partial X_2}{\partial \theta} & 
\fraco{1}{2}\ten{S}_{12}\fraco{\partial X_3}{\partial \lambda} + \fraco{1}{2}\ten{S}_{21}\fraco{\partial X_3}{\partial \lambda} + \ten{S}_{22}\fraco{\partial X_3}{\partial \theta} \\[8 pt]
\end{pmatrix}.
\end{align}
By the symmetry of $\ten{S}$, we have
\begin{equation}
\ten{S} \cddot \frac{\partial \gamma}{\partial (\mathbf{\nabla_{\theta}X})^\textrm{T}} = \begin{pmatrix}
\ten{S}_{11}\fraco{\partial X_1}{\partial \lambda} + \ten{S}_{12}\fraco{\partial X_1}{\partial \theta} & \ten{S}_{11}\fraco{\partial X_2}{\partial \lambda} + \ten{S}_{12}\fraco{\partial X_2}{\partial \theta} & \ten{S}_{11}\fraco{\partial X_3}{\partial \lambda} + \ten{S}_{12}\fraco{\partial X_3}{\partial \theta} \\[10 pt]
\ten{S}_{21}\fraco{\partial X_1}{\partial \lambda} + \ten{S}_{22}\fraco{\partial X_1}{\partial \theta} & 
\ten{S}_{21}\fraco{\partial X_2}{\partial \lambda} + \ten{S}_{22}\fraco{\partial X_2}{\partial \theta} & 
\ten{S}_{21}\fraco{\partial X_3}{\partial \lambda} + \ten{S}_{22}\fraco{\partial X_3}{\partial \theta} \\
\end{pmatrix}.
\end{equation}
From the dot product,
\begin{align}
\ten{S} \cdot (\mathbf{\nabla_{\theta}X})^\textrm{T} & = \left(\ten{S}_{11} \vec{e_1}\vec{e_1} + \ten{S}_{12}\vec{e_1}\vec{e_2} + \ten{S}_{21}\vec{e_2}\vec{e_1} + \ten{S}_{22}\vec{e_2}\vec{e_2} \right) \cdot \nonumber \\[6 pt]
&  \left(\frac{\partial X_1}{\partial \lambda} \vec{e_1}\vec{e_1} + \frac{\partial X_2}{\partial \lambda} \vec{e_1}\vec{e_2} + \frac{\partial X_3}{\partial \lambda} \vec{e_1}\vec{e_3} + \frac{\partial X_1}{\partial \theta} \vec{e_2}\vec{e_1} + \frac{\partial X_2}{\partial \theta} \vec{e_2}\vec{e_2} + \frac{\partial X_3}{\partial \theta} \vec{e_2}\vec{e_3}\right)\\[6 pt]
& = \left(\begin{matrix}
\ten{S}_{11}\fraco{\partial X_1}{\partial \lambda} + \ten{S}_{12}\fraco{\partial X_1}{\partial \theta} & \ten{S}_{11}\fraco{\partial X_2}{\partial \lambda} + \ten{S}_{12}\fraco{\partial X_2}{\partial \theta} & \ten{S}_{11}\fraco{\partial X_3}{\partial \lambda} + \ten{S}_{12}\fraco{\partial X_3}{\partial \theta} \\[10 pt]
\ten{S}_{21}\fraco{\partial X_1}{\partial \lambda} + \ten{S}_{22}\fraco{\partial X_1}{\partial \theta} & 
\ten{S}_{21}\fraco{\partial X_2}{\partial \lambda} + \ten{S}_{22}\fraco{\partial X_2}{\partial \theta} & 
\ten{S}_{21}\fraco{\partial X_3}{\partial \lambda} + \ten{S}_{22}\fraco{\partial X_3}{\partial \theta} \\
\end{matrix}\right) \\[6 pt]
& = \ten{S} \cddot \frac{\partial \gamma}{\partial (\mathbf{\nabla_{\theta}X})^\textrm{T}},
\end{align}
so that we recover Eq.\ \eqref{eq:tenchain}. 
\subsection{Chain rule expansion of $\ten{S}$}
As observed in Eq.\ \eqref{eq:sexpansion}, the second Piola-Kirchhoff stress tensor $\ten{S}$ can be computed by observing that $\C = \ten{I}_2 + 2\gamma \Gz^{-1}$, where $\ten{I}_2$ is the rank 2 identity tensor. Applying the chain rule twice yields
\begin{equation}
\ten{S} = \frac{\partial W}{\partial \gamma} = \frac{\partial W}{\partial \C} \cddot \frac{\partial \C}{\partial \gamma} = 2\frac{\partial W}{\partial \C} \cdot \Gz^{-1},
\end{equation}
and
\begin{equation}
\frac{\partial W}{\partial \C} = \frac{\partial W}{\partial I_1}\frac{\partial I_1}{\partial \C} + \frac{\partial W}{\partial I_2}\frac{\partial I_2}{\partial \C}.
\end{equation}The rank 2 tensor $\fraco{\partial I_j}{\partial \C}$ is obtained by taking the derivative of the invariant $I_j$ with respect to each element of $\C$. Carrying out these calculations, 
\begin{gather}
I_1  = \Tr{\C} - 2 = \C_{11}+\C_{22} - 2, \\[6 pt]
\begin{matrix}
\fraco{\partial{I_1}}{\partial \C_{11}}  = 1, & \fraco{\partial{I_1}}{\partial \C_{12}} = 0,\\[8 pt]
\fraco{\partial{I_1}}{\partial \C_{21}} = 0, & \fraco{\partial{I_1}}{\partial \C_{22}} = 1,
\end{matrix}\\[8 pt]
\fraco{\partial I_1}{\partial \C} = \ten{I}_2,
\end{gather}
\begin{gather}
I_2  = \Det{\C} - 1 = \C_{11}\C_{22} - \C_{12}\C_{21} - 1, \\[6 pt]
\begin{matrix}
\fraco{\partial{I_2}}{\partial \C_{11}}  = \C_{22}, & \fraco{\partial{I_2}}{\partial \C_{12}} = -\C_{21},\\[8 pt]
\fraco{\partial{I_2}}{\partial \C_{21}} = -\C_{12}, & \fraco{\partial{I_2}}{\partial \C_{22}} = \C_{11},
\end{matrix}\\[8 pt]
\frac{\partial{I_2}}{\partial \C} = (\Det{\C})(\C^{-1})^\textrm{T} = (\Det{\C})\left((\G \Gz^{-1})^{-1}\right)^\textrm{T} = (\Det{\C})\left(\Gz \G^{-1}\right)^\textrm{T} = (\Det{\C})\left(\G^{-1}\Gz\right).
\end{gather}
The last equality follows from the symmetry of $\G$ and $\Gz$. 
At last we have that 
\begin{align}
\ten{S}  & = 2\frac{\partial W}{\partial \C} \cdot \Gz^{-1} = 2\left(\frac{\partial W}{\partial I_1}\frac{\partial I_1}{\partial \C} + \frac{\partial W}{\partial I_2}\frac{\partial I_2}{\partial \C} \right) \cdot \Gz^{-1}\\[6 pt]
& =  2\left(\frac{\partial W}{\partial I_1}\ten{I}_2 + \frac{\partial W}{\partial I_2}(\Det{\C})\left(\G^{-1}\Gz\right) \right) \cdot \Gz^{-1}\\[6 pt] 
\ten{S} & = 2\frac{\partial W}{\partial I_1}\Gz^{-1} + 2\frac{\partial W}{\partial I_2}(\Det{\C})\G^{-1}.
\end{align}
\subsection{The force density}
Computing the force density in Eq.\ \eqref{eq:feqdivdotPt}, we use the chain rule expansions of $\ten{P}^{\textrm{T}}$ and $\ten{S}$:
\begin{align}
\nabla_\theta \, \cdot \, & = \frac{1}{\sqrt{\Det{\Gz}}} \left(\frac{\partial}{\partial \lambda} \sqrt{\Det{\Gz}} \textrm{ } \vec{e_1} + \frac{\partial}{\partial \theta} \sqrt{\Det{\Gz}} \textrm{ } \vec{e_2}\right) \cdot ,\\[6 pt]
\vec{F} & =\nabla_\theta \cdot \ten{P}^{\textrm{T}} = \nabla_\theta \cdot \left(\ten{S} \cdot (\mathbf{\nabla_{\theta}X})^{\textrm{T}}\right),\\[6 pt]\nonumber
\vec{F}& = \nabla_\theta \cdot 
\Bigg(\left(\ten{S}_{11}\frac{\partial X_1}{\partial \lambda} + \ten{S}_{12}\frac{\partial X_1}{\partial \theta} \right)\vec{e_1} \vec{e_1} + 
\left(\ten{S}_{11}\frac{\partial X_2}{\partial \lambda} + \ten{S}_{12}\frac{\partial X_2}{\partial \theta} \right)\vec{e_1} \vec{e_2}\\[6 pt] 
& + \left(\ten{S}_{11}\frac{\partial X_3}{\partial \lambda} + \ten{S}_{12}\frac{\partial X_3}{\partial \theta} \right)\vec{e_1} \vec{e_3} + \left(\ten{S}_{21}\frac{\partial X_1}{\partial \lambda} + \ten{S}_{22}\frac{\partial X_1}{\partial \theta} \right)\vec{e_2} \vec{e_1} \\[6 pt]\nonumber
& + \left(\ten{S}_{21}\frac{\partial X_2}{\partial \lambda} + \ten{S}_{22}\frac{\partial X_2}{\partial \theta} \right)\vec{e}_2 \vec{e}_2 + 
\left(\ten{S}_{21}\frac{\partial X_3}{\partial \lambda} + \ten{S}_{22}\frac{\partial X_3}{\partial \theta}\right)\vec{e}_2 \vec{e}_3\Bigg).
\end{align}
Thus
\begin{align}
\vec{F} & = \frac{1}{\sqrt{\Det{\Gz}}} \begin{pmatrix}
\fraco{\partial}{\partial \lambda} \sqrt{\Det{\Gz}} \left(\ten{S}_{11}\fraco{\partial X_1}{\partial \lambda} + \ten{S}_{12}\fraco{\partial X_1}{\partial \theta} \right) + \fraco{\partial}{\partial \theta} \sqrt{\Det{\Gz}}\left(\ten{S}_{21}\fraco{\partial X_1}{\partial \lambda} + \ten{S}_{22}\frac{\partial X_1}{\partial \theta} \right) \\[10 pt]
 \fraco{\partial}{\partial \lambda} \sqrt{\Det{\Gz}} \left(\ten{S}_{11}\fraco{\partial X_2}{\partial \lambda} + \ten{S}_{12}\fraco{\partial X_2}{\partial \theta} \right) + \fraco{\partial}{\partial \theta} \sqrt{\Det{\Gz}}\left(\ten{S}_{21}\fraco{\partial X_2}{\partial \lambda} + \ten{S}_{22}\fraco{\partial X_2}{\partial \theta} \right) \\[10 pt]
 \fraco{\partial}{\partial \lambda} \sqrt{\Det{\Gz}} \left(\ten{S}_{11}\fraco{\partial X_3}{\partial \lambda} + \ten{S}_{12}\fraco{\partial X_3}{\partial \theta} \right) + \fraco{\partial}{\partial \theta} \sqrt{\Det{\Gz}}\left(\ten{S}_{21}\fraco{\partial X_3}{\partial \lambda} + \ten{S}_{22}\fraco{\partial X_3}{\partial \theta} \right)
\end{pmatrix}\\[12 pt]
\label{eq:finalF}
\vec{F} & =  \frac{1}{\sqrt{\Det{\Gz}}} \left(\frac{\partial}{\partial \lambda} \sqrt{\Det{\Gz}} \left(\ten{S}_{11}\frac{\partial \vec{X}}{\partial \lambda} + \ten{S}_{12}\frac{\partial \vec{X}}{\partial \theta} \right) + \frac{\partial}{\partial \theta} \sqrt{\Det{\Gz}}\left(\ten{S}_{21}\frac{\partial \vec{X}}{\partial \lambda} + \ten{S}_{22}\frac{\partial \vec{X}}{\partial \theta} \right)\right).
\end{align}
Eq.\ \eqref{eq:finalF} recovers Eq.\ \eqref{eq:vdforce} from the main text.
%
%
%
\section*{References}
\bibliographystyle{elsarticle-num}
\bibliography{WandaBib}
\end{document}